\documentclass[preprint,numbers,sort&compress]{elsarticle}
\usepackage{amsmath,amsfonts,amsthm,mathabx}
\usepackage{graphicx}
\usepackage{upquote}
\usepackage{pgfplots}
\pgfplotsset{compat=newest}
\usetikzlibrary{plotmarks}
\usetikzlibrary{arrows.meta}
\usetikzlibrary{spy}
\usepgfplotslibrary{patchplots}
\usepackage{grffile}
\usepackage{subfloat}
\newcommand{\bx}{{\boldsymbol x}}
\newcommand{\bu}{{\boldsymbol u}}
\newcommand{\bv}{{\boldsymbol v}}
\newcommand{\bw}{{\boldsymbol w}}
\newcommand{\bg}{{\boldsymbol g}}
\newcommand{\bU}{{\boldsymbol U}}
\newcommand{\bV}{{\boldsymbol V}}
\newcommand{\bG}{{\boldsymbol G}}
\newcommand{\bW}{{\boldsymbol W}}
\newcommand{\rme}{\mathrm{e}}

\newcommand{\RR}{\mathbb{R}}
\newcommand{\CC}{\mathbb{C}}

\date{}

\definecolor{rev1}{RGB}{255,0,0} 
\definecolor{rev3}{RGB}{0,0,255} 
\definecolor{rev4}{RGB}{0,153,0} 
\definecolor{our}{RGB}{255,128,0} 

\usepackage{hyperref}
\begin{document}
\begin{frontmatter}
\title{A second order directional split exponential integrator for
       systems of advection--diffusion--reaction equations}

\author[label1]{Marco Caliari\corref{cor1}}
\ead{marco.caliari@univr.it}
\cortext[cor1]{Corresponding author}
\affiliation[label1]{organization={Department of Computer Science, University
    of Verona},
  country={Italy}}
\myfooter[L]{} 
\author[label1]{Fabio Cassini}
\ead{fabio.cassini@univr.it}
\begin{abstract}
  We propose a second order exponential scheme suitable for two-component
  coupled systems of stiff {\color{black} evolutionary}
  advection--diffusion--reaction equations in two and
  three space dimensions. It is based on a directional splitting of the
  involved matrix functions, which allows for a simple yet efficient
  implementation
  through the computation of small-sized exponential-like
  functions and tensor-matrix products.
  The procedure straightforwardly extends to the case of an arbitrary number
  of components and to any space dimension. Several numerical examples
  in 2D and 3D with physically relevant (advective)
  Schnakenberg, FitzHugh--Nagumo, DIB,
  and advective Brusselator models
  clearly show the advantage of the approach against
  state-of-the-art techniques.
\end{abstract}

\end{frontmatter}
\section{Introduction}
Efficiently solving systems of Advection--Diffusion--Reaction (ADR) equations
is of paramount importance in practical applications. In fact, these kinds of
models effectively capture many physical, chemical, and biological phenomena,
such as
electrodeposition~\cite{BLS13}, biochemical reactions~\cite{MWM03,TGC99},
electric current flows~\cite{GLRS19}, {\color{black}tumor growth~\cite{CS01}, and
epidemic dynamics~\cite{MPV08},} among the others. From a numerical
point of view, the computation of approximated solutions gives rise to many
challenges. For example, in this context there is often the intrinsic need
of a fine spatial grid in order to correctly capture the dynamics of
the system (the formation of the so-called Turing patterns~\cite{DASS20}, for
instance). This,
in turn, typically translates into working in a \emph{stiff} regime~\cite{HW96},
which requires a careful choice of the underlying numerical schemes.

In particular, in this work we focus on the numerical integration of
\emph{two-component} systems of ADR equations in the form
\begin{equation}\label{eq:pdes}
  \left\{
  \begin{aligned}
    \partial_t u(t,\bx)&=\mathcal{K}^uu(t,\bx)+g^u(u(t,\bx),v(t,\bx)),\\
    \partial_t v(t,\bx)&=\mathcal{K}^vv(t,\bx)+g^v(u(t,\bx),v(t,\bx)),
  \end{aligned}\right.
\end{equation}
although the method presented later on in the manuscript easily extends to an
arbitrary number of components.
Here $u,v\colon[0,T]\times\Omega\subset {\color{black}\RR\times\RR^d}\to\RR$
represent the unknowns,
$\mathcal{K}^u,\mathcal{K}^v$ are linear
advection--diffusion operators, while
$g^u, g^v$ are the nonlinear reaction terms. The choice of the latter basically
determines the model and phenomenon under investigation.
We assume that the spatial domain $\Omega$ is the Cartesian product
of one-dimensional
intervals, that is $\Omega=[a_1,b_1]\times\cdots\times [a_d,b_d]${\color{black},
as it is commonly considered in ADR examples from the literature (see,
for instance, References~\cite{AMS23,BTMSC23,GLRS19,DASS20,BKW18}).
The system is finally completed with appropriate initial conditions
and with homogeneous Neumann boundary
conditions, which appear to be widely used by the scientific community.}
We introduce a spatial grid of size $n_1\times\cdots \times n_d$ and
apply the Method Of Lines (MOL) to
model~\eqref{eq:pdes} {\color{black}to get a system with \textit{Kronecker
sum} structure. In fact, we suppose to obtain}
a system of stiff ordinary differential equations in the form
\begin{equation}\label{eq:odes}
  \left\{\begin{aligned}
  \bu'(t)&=K^{\bu}\bu(t)+\bg^{\bu}(\bu(t),\bv(t)),\\
  \bv'(t)&=K^{\bv}\bv(t)+\bg^{\bv}(\bu(t),\bv(t)),
  \end{aligned}\right.
\end{equation}
where $K^{\bu}$ and $K^{\bv}$ are
matrices of size $N\times N$, with $N=n_1\cdots n_d$, that discretize the
linear operators $\mathcal{K}^u$ and $\mathcal{K}^v$, respectively.
The matrices $K^{\bu}$ and $K^{\bv}$ are Kronecker sums.
By definition, a matrix $K\in\CC^{N\times N}$ is a Kronecker sum if it
can be decomposed as
\begin{subequations}\label{eq:kronsum}
\begin{equation}\label{eq:kronsumK}
  K=A_d\oplus A_{d-1}\oplus \cdots \oplus A_1=
  \sum_{\mu=1}^d A_{\otimes\mu},
\end{equation}
where
\begin{equation}\label{eq:kronAmu}
  A_{\otimes \mu}=I_d\otimes \cdots\otimes I_{\mu+1}\otimes
  A_\mu\otimes I_{\mu-1}\otimes \cdots\otimes I_1.
\end{equation}
\end{subequations}
{\color{black}
In our context $I_\mu$ and $A_\mu$, for $\mu=1,\ldots,d$,
are matrices of \emph{small} size $n_\mu\times n_\mu$
and represent the identity and a one-dimensional 
linear differential operator along the direction $\mu$, respectively.}
The symbol $\otimes$ denotes the standard Kronecker product between matrices.
{\color{black}The Kronecker sum structure
  typically arises when applying finite
differences discretizations~\cite{CCEOZ22,CC23} to the differential
operators involved in the equations, but for instance 
tensor product finite elements~\cite{MMPC22,CMM23} fit into this class as
well.}
{\color{black}Notice that, for the method we are going to propose,
  we do not assume the matrices
  $A_\mu$ having a specific sparsity pattern.}

Several methods can be employed for the time
integration of system~\eqref{eq:odes}.
{\color{black}In presence of} stiffness,
an attractive choice are integrators of exponential or implicit type,
that are useful
to alleviate the stability restrictions on the time step
size that come with explicit methods.
Of course, ad-hoc techniques have to be adopted in order to
retain efficiency of computation.
For instance, among schemes of exponential type~\cite{HO10},
Lawson integrators~\cite{L67}
(also known as Integration Factor methods~\cite{JZ16}) are appealing,
since they require
just the action of the matrix exponential
function and thus they can directly exploit the Kronecker sum structure of the
system (see Reference~\cite{CC23}, and Reference~\cite{CCEOZ22} for
some insights on GPU scalability). On the other hand, it is well-known
that they may suffer from a bad convergence behavior, in particular when the
boundary conditions are
non-periodic (see Reference~\cite{HLO20}).
Different exponential integrators rely on general purpose
Krylov methods for the
approximation
of the matrix exponential or related functions with a fixed size of the Krylov
space (see Reference~\cite{BKW18}) or with incomplete orthogonalization
(see References~\cite{LPR19,GRT18}). More recent techniques that exploit the
underlying Kronecker sum structure are presented, for instance, in
References~\cite{AAKW20,CCZ23phi,CMM23,MMPC22,DASS20}.
Concerning implicit methods, in the two-dimensional case
low-rank or other Model
Order Reduction techniques combined
with a matrix formulation
of the problem can be applied
to classical schemes such as IMplicit EXplicit
methods (see References~\cite{DASS20,AMS23}).
We finally mention also the possibility
to employ explicit schemes such as
a Total Variation Diminishing explicit Runge--Kutta
scheme (see Reference~\cite{SMTT23}), provided that a sufficiently
small time step size is taken.

In this manuscript, we propose an exponential scheme
of second order for system~\eqref{eq:odes}.
It is a \emph{directional split} version of the well-known exponential
Runge--Kutta method of order two ETD2RK, and it approximates
the required four actions of matrix $\varphi_\ell$ functions
by means of four tensor-matrix products (using the so-called Tucker operator).
This approach allows for a simple yet very efficient implementation, thanks to
the high performance level 3 BLAS available for basically any kind of modern
computer system. {\color{black}A preliminary study of the
  effectiveness of the directional
splitting of the matrix $\varphi_\ell$ functions in exponential integrators
has been done in Reference~\cite{CC23}, but there the numerical examples
were limited to recover the steady state of an algebraic Riccati
equation and to measure its performance on a prototypical evolutionary
equation (whose explicit analytical solution was known).
Here, the focus is to show that the proposed exponential
scheme is highly competitive with state-of-the-art methods in
the solution of important ADR models, and in particular for the
challenging task of recovering the arising patterns.}

The remaining part of paper is organized as follows.
In Section~\ref{sec:expint} we will
briefly recall integrators of exponential type
and the tensor framework needed to implement them efficiently.
In Section~\ref{sec:method} we will present ETD2RKds, which is the proposed
integrator with directional splitting.
{\color{black}Although the method is presented in arbitrary dimension $d$,}
we proceed in Section~\ref{sec:numexp} by testing it on a wide variety of
two-dimensional and three-dimensional numerical examples
with popular models in the literature, i.e.,
the {\color{black}(advective) Schnakenberg, the FitzHugh--Nagumo, the DIB,
  and the advective
Brusselator} models.
{\color{black}We will see that the proposed approach
is able to effectively and efficiently integrate numerically the problems,
both in terms of achieved accuracy and wall-clock time, in comparison with
other popular techniques. Moreover, we will
also show that we can attain the expected
patterns of the underlying models in a short amount of time.}
Finally, we will draw the conclusions in Section~\ref{sec:conc}.

\section{Exponential integrators for
  systems of partial differential equations}\label{sec:expint}
The starting point for the exponential integrators under consideration is the
variation-of-constants formula applied to system~\eqref{eq:odes} with time step
size $\tau=t_{n+1}-t_n$, i.e.,
\begin{equation}\label{eq:voc}
\begin{aligned}
  \bu(t_{n+1})&=\rme^{\tau K^{\bu}}\bu(t_n)+\tau\int_0^1
  \rme^{(1-\theta)\tau K^{\bu}}\bg^{\bu}(\bu(t_n+\theta\tau),\bv(t_n+\theta\tau))d\theta,\\
  \bv(t_{n+1})&=\rme^{\tau K^{\bv}}\bv(t_n)+\tau\int_0^1
  \rme^{(1-\theta)\tau K^{\bv}}\bg^{\bv}(\bu(t_n+\theta\tau),\bv(t_n+\theta\tau))d\theta.\\
\end{aligned}
\end{equation}
We then apply the trapezoidal quadrature rule, by introducing an explicit
intermediate
stage to approximate the evaluation at $t_{n+1}$
(see also~Reference~\cite{CC23}).
Putting everything together, we get the second order Lawson scheme
\begin{equation*}
\begin{aligned}
  \bu_{n2}&=\rme^{\tau K^{\bu}}(\bu_n+\tau\bg^{\bu}(\bu_n,\bv_n)),\\
  \bv_{n2}&=\rme^{\tau K^{\bv}}(\bv_n+\tau\bg^{\bv}(\bu_n,\bv_n)),\\
  \bu_{n+1}&=\rme^{\tau K^{\bu}}\left(\bu_n+\frac{\tau}{2}\bg^{\bu}(\bu_n,\bv_n)\right)+
  \frac{\tau}{2}\bg^{\bu}(\bu_{n2},\bv_{n2}),\\
  \bv_{n+1}&=\rme^{\tau K^{\bv}}\left(\bv_n+\frac{\tau}{2}\bg^{\bv}(\bu_n,\bv_n)\right)+
  \frac{\tau}{2}\bg^{\bv}(\bu_{n2},\bv_{n2}).
\end{aligned}
\end{equation*}
The same method can be recovered also by multiplying the two
equations of system~\eqref{eq:odes} by the integration
factors $\rme^{-t K^\bu}$ and $\rme^{-t K^\bv}$, respectively, and applying the
explicit Heun method (see for instance Reference~\cite{JZ16}).
The just mentioned Lawson scheme can be efficiently implemented by
directly exploiting the
Kronecker
sum structure of the involved matrices. Indeed,
{\color{black}in the two-dimensional case (that is,
  for a matrix $\tau K=\tau A_2\oplus \tau A_1$), we can use
  the well-known formula (see, for instance, Reference~\cite{N69})
\begin{equation*}
\rme^{\tau K}\bw=\mathrm{vec}\left(\rme^{\tau A_1}\bW\left(\rme^{\tau A_2}\right)^{\sf T}\right),
\end{equation*}
where $\bW$ is a matrix of size $n_1\times n_2$ and
$\mathrm{vec}$ is
the operator which stacks the columns of the input matrix into a suitable
column vector such that $\mathrm{vec}(\bW)=\bw$. This allows for computing
the action of the matrix exponential without assembling the matrix $K$
and by performing two products with the dense
matrices $\rme^{\tau A_1}$ and $\rme^{\tau A_2}$, respectively.}
For a matrix $\tau K$ with $d$-dimensional
Kronecker sum structure~\eqref{eq:kronsum},
the formula above can be generalized
by employing
the representation
\begin{equation}\label{eq:Tucker}
  \rme^{\tau K}\bw = \mathrm{vec}\left(\bW
  \times_1\rme^{\tau A_1}\times_2\cdots
  \times_d\rme^{\tau A_d}\right)=
  \mathrm{vec}\left(\bW\bigtimes_{\mu=1}^d\rme^{\tau A_\mu}\right),
\end{equation}
{\color{black}see Reference~\cite{CCEOZ22}}, where $\bW$
is an order-$d$ tensor of size $n_1\times\cdots\times n_d$.
Here, $\times_\mu$ denotes the \emph{$\mu$-mode product}, {\color{black}
  which  multiplies the matrix $\rme^{\tau A_\mu}$ onto the $\mu$-fibers
  (i.e., the generalizations to tensors of
  matrix columns and rows) of
  the tensor~$\bW$.}
{\color{black}Notice that also in the $d$-dimensional case it can be
efficiently realized  by means of high performance
level 3 BLAS operations at
a computational cost $\mathcal{O}(n_\mu N)$,
since the matrix
  exponentials $\rme^{\tau A_\mu}$ are dense matrices.}
The concatenation of $d$ different $\mu$-mode products,
that we denoted as $\bigtimes_{\mu=1}^d$, is called {\color{black}a}
\emph{Tucker operator},
and clearly has computational cost $\mathcal{O}((n_1+\cdots+n_d) N)$.
The exponentials of the matrices $\tau A_\mu$
 can be efficiently computed by polynomial  or rational
 approximations (see~References~\cite{AMH10,SID19,CZ19,AIDAJ23})
 together
 with a scaling and squaring algorithm,
  or even by a spectral decomposition (see~Reference~\cite{DASS20}).
We refer the reader also to Reference~\cite{CCZ23kp} for more details about
formula~\eqref{eq:Tucker}, the Tucker operator and related operations
(implemented in a MATLAB package\footnote{Available
  at \url{https://github.com/caliarim/KronPACK},
  commit \texttt{562a9da}.}). 

{\color{black}Using the tensor formalism just introduced, it is also possible
  to equivalently rewrite the vector ODE system~\eqref{eq:odes} as
\begin{equation}\label{eq:todes}
  \left\{\begin{aligned}
  \bU'(t)&=\sum_{\mu=1}^d (\bU(t)\times_\mu A_\mu^{\bu})+
  \bG^{\bu}(\bU(t),\bV(t)),\\
  \bV'(t)&=\sum_{\mu=1}^d (\bV(t)\times_\mu A_\mu^{\bv})+
  \bG^{\bv}(\bU(t),\bV(t)),
  \end{aligned}\right.
\end{equation}
where
\begin{align*}
  \mathrm{vec}(\bU(t))&=\bu(t),&
  \mathrm{vec}(\bG^{\bu}(\bU(t),\bV(t)))&=\bg^{\bu}(\bu(t),\bv(t)),\\
  \mathrm{vec}(\bV(t))&=\bv(t),&
  \mathrm{vec}(\bG^{\bv}(\bU(t),\bV(t)))&=\bg^{\bv}(\bu(t),\bv(t)).
\end{align*}
Notice that, in the formulation above, the actions of
$K^\bu$ and $K^\bv$ are expressed in tensor form
thanks to the identities
\begin{equation*}
  K^\bu \bu(t)=\mathrm{vec}\left(\sum_{\mu=1}^d \bU(t)\times_\mu A^\bu_\mu\right)
  \quad \text{and}
  \quad
    K^\bv \bv(t)=\mathrm{vec}\left(\sum_{\mu=1}^d \bV(t)\times_\mu A^\bv_\mu\right).
\end{equation*}
Each of these two evaluations has a computational
cost equivalent to that of a single
Tucker operator.
Then, the aforementioned second order Lawson scheme
in tensor formulation is
\begin{equation}\label{eq:law2bten}
\begin{aligned}
  \bU_{n2}&=
  \left(\bU_n+\tau\bG^{\bu}(\bU_n,\bV_n)\right)\bigtimes_{\mu=1}^d\rme^{\tau A^\bu_\mu},\\
  \bV_{n2}&=\left(\bV_n+\tau\bG^{\bv}(\bU_n,\bV_n)\right)\bigtimes_{\mu=1}^d\rme^{\tau A^\bv_\mu},\\
  \bU_{n+1}&=\left(\bU_n+\frac{\tau}{2}\bG^{\bu}(\bU_n,\bV_n)\right)\bigtimes_{\mu=1}^d\rme^{\tau A^\bu_\mu}+
  \frac{\tau}{2}\bG^{\bu}(\bU_{n2},\bV_{n2}),\\
  \bV_{n+1}&=\left(\bV_n+\frac{\tau}{2}\bG^{\bv}(\bU_n,\bV_n)\right)\bigtimes_{\mu=1}^d\rme^{\tau A^\bv_\mu}+
  \frac{\tau}{2}\bG^{\bv}(\bU_{n2},\bV_{n2}).
\end{aligned}
\end{equation}
For the two-dimensional case, the scheme can be written without
introducing tensor notation as
    \begin{equation*}
\begin{aligned}
  \bU_{n2}&=\rme^{\tau A^\bu_1}\left(\bU_n+\tau   \bG^\bu(\bU_n,\bV_n)\right)
  \left(\rme^{\tau A_2^\bu}\right)^{\sf T},\\
  \bV_{n2}&=\rme^{\tau A^\bv_1}\left(\bV_n+\tau   \bG^\bv(\bU_n,\bV_n)\right)
  \left(\rme^{\tau A_2^\bv}\right)^{\sf T},\\
    \bU_{n+1}&=\rme^{\tau A^\bu_1}\left(\bU_n+
    \frac{\tau}{2}   \bG^\bu(\bU_n,\bV_n)\right)
  \left(\rme^{\tau A_2^\bu}\right)^{\sf T}+\frac{\tau}{2}\bG^\bu(\bU_{n2},\bV_{n2}),\\
    \bV_{n+1}&=\rme^{\tau A^\bv_1}\left(\bV_n+
    \frac{\tau}{2}   \bG^\bv(\bU_n,\bV_n)\right)
  \left(\rme^{\tau A_2^\bv}\right)^{\sf T}+\frac{\tau}{2}\bG^\bv(\bU_{n2},\bV_{n2}).\\
\end{aligned}
\end{equation*}
We will refer to method~\eqref{eq:law2bten} as Lawson2b.}

Different approximations in the integrals of formula~\eqref{eq:voc}
lead to other exponential integrators.
In particular, if we interpolate {\color{black}at $\theta=0$ and
  $\theta=1$ the nonlinear functions
  $\bg^{\bu}$ and $\bg^{\bv}$ by a polynomial,
and we introduce an intermediate
stage for the approximation at $\theta=1$, we get the
exponential Runge--Kutta scheme of order two
\begin{equation}\label{eq:erk}
\begin{aligned}
  \bu_{n2}&=\bu_n+\tau\varphi_1(\tau K^{\bu})(K^{\bu}\bu_n+\bg^{\bu}(\bu_n,\bv_n)),\\
  \bv_{n2}&=\bv_n+\tau\varphi_1(\tau K^{\bv})(K^{\bv}\bv_n+\bg^{\bv}(\bu_n,\bv_n)),\\
  \bu_{n+1}&=\bu_{n2}+\tau\varphi_2(\tau K^{\bu})
  \left(\bg^{\bu}(\bu_{n2},\bv_{n2})-\bg^{\bu}(\bu_n,\bv_n)\right),\\
  \bv_{n+1}&=\bv_{n2}+\tau\varphi_2(\tau K^{\bv})
  \left(\bg^{\bv}(\bu_{n2},\bv_{n2})-\bg^{\bv}(\bu_n,\bv_n)\right),
\end{aligned}
\end{equation}
known as ETD2RK (see, for instance, Reference~\cite{BSW05}).}
The matrix $\varphi_\ell$ functions, for a generic
matrix $X\in\CC^{N\times N}$, are defined by
\begin{equation}\label{eq:phiell}
  \varphi_\ell(X)=\int_0^1\frac{\theta^{\ell-1}}{(\ell-1)!}\rme^{(1-\theta)X}d\theta,
  \quad \ell>0.
\end{equation}
Several techniques are available in the literature to compute the matrix
$\varphi_\ell$ functions (see References~\cite{BSW05,LYL22,CC23} for methods
based on rational or polynomial approximations or quadrature formulas,
{\color{black}suitable for small sized matrices})
or their action to vectors (see References~\cite{CCZ23,GRT18,LPR19} for
Krylov methods with incomplete orthogonalization, {\color{black}suitable
  for large sized and sparse matrices}).

Unfortunately, even if $X=\tau K$, property~\eqref{eq:Tucker} is not valid
for the matrix $\varphi_\ell$ functions.
However, it is still possible {\color{black}to} exploit the Kronecker sum structure of $K$ to
compute the action
of the matrix $\varphi_\ell$ function on a vector $\bw$ efficiently.
In fact, if we introduce a quadrature
formula for integral definition~\eqref{eq:phiell}, we can approximate
the desired
quantity by computing the actions $\rme^{(1-\theta_i)\tau K}\bw$ using Tucker
operators, being $\theta_i$
the quadrature nodes. In this way, the approximation
via quadrature formula requires one Tucker operator for each quadrature node.
The total number of quadrature nodes can be kept reasonably small by employing
a modified scaling and squaring
technique for the $\varphi_\ell$ functions.
Each iteration of the squaring procedure requires then the action of a matrix
exponential, which can be realized again by a Tucker operator.
{\color{black}Such an approach is employed in
  Reference~\cite{CCZ23phi} and is called
\textsc{phiks} (PHI-functions of Kronecker Sums), while a similar
  technique is presented in Reference~\cite{CMM23}.}
  We refer to those manuscripts for the
  details on the choice of the quadrature formulas and 
  of the scaling parameter.

\section{Directional splitting of the matrix
  \texorpdfstring{$\varphi_\ell$}{phi} functions}\label{sec:method}
As already mentioned, the splitting property
\begin{equation*}
  \rme^{\tau K}\bw = \rme^{\tau A_{\otimes 1}}\cdots\rme^{\tau A_{\otimes d}}\bw,
\end{equation*}
which is at the basis of equivalence~\eqref{eq:Tucker}, does not extend
directly to the $\varphi_\ell$ functions.
However, it has been observed in Reference~\cite{CC23}
that a simple directional
splitting
of the matrix $\varphi_\ell$ functions, that is
\begin{subequations}\label{eq:phisplitcomp}
\begin{equation}\label{eq:phisplit}
  \varphi_\ell(\tau K)\bw=
  (\ell!)^{d-1}\varphi_\ell(\tau A_{\otimes 1})
  \cdots\varphi_\ell(\tau A_{\otimes d})\bw+
  \mathcal{O}(\tau^2),
\end{equation}
  gives an approximation
  compatible with second order integrators. The advantage of using the above
  approximation is that it can be efficiently computed in tensor form
  (similarly to the matrix exponential case) as
\begin{equation}\label{eq:phisplit2}
  \varphi_\ell(\tau A_{\otimes 1})
  \cdots\varphi_\ell(\tau A_{\otimes d})\bw=
  \mathrm{vec}\left(\bW\bigtimes_{\mu=1}^d\varphi_\ell(\tau A_\mu)\right).
\end{equation}
\end{subequations}
The approximation introduced by the directional
splitting can affect in general the magnitude of the
error when applied in exponential integrators.
However, we will see later in the numerical examples that the overall
  computational cost required to reach the same level of accuracy is smaller.
  This is possible since only small sized matrix $\varphi_\ell$ functions
  have to be computed, and the Tucker operator
  (which can be efficiently realized as explained in Section~\ref{sec:expint})
  has to be applied.

In particular, if we embed splitting formulas~\eqref{eq:phisplitcomp} in
the ETD2RK integrator~\eqref{eq:erk},
we obtain a scheme that we call \emph{directional split ETD2RK}, and more
concisely ETD2RKds.
  The complete expression of the proposed scheme
  for the solution of the tensor ODE system~\eqref{eq:todes} is then
\begin{equation}\label{eq:etd2rk_d_split}
\begin{aligned}
  \bU_{n2}&=\bU_n+\tau
  \left(\sum_{\mu=1}^d(\bU_n\times_\mu A_\mu^\bu)+
  \bG^\bu(\bU_n,\bV_n)\right)
  \bigtimes_{\mu=1}^d\varphi_1(\tau A^\bu_\mu),\\
  \bV_{n2}&=\bV_n+\tau
  \left(\sum_{\mu=1}^d(\bV_n\times_\mu A_\mu^\bv)+
  \bG^\bv(\bU_n,\bV_n)\right)
\bigtimes_{\mu=1}^d\varphi_1(\tau A^\bv_\mu),\\
  \bU_{n+1}&=\bU_{n2}+2^{d-1}\tau (\bG^\bu(\bU_{n2},\bV_{n2})
  -\bG^\bu(\bU_n,\bV_n))
\bigtimes_{\mu=1}^d\varphi_2(\tau A^\bu_\mu),\\
  \bV_{n+1}&=\bV_{n2}+2^{d-1}\tau (\bG^\bv(\bU_{n2},\bV_{n2})
  -\bG^\bv(\bU_n,\bV_n))
\bigtimes_{\mu=1}^d\varphi_2(\tau A^\bv_\mu),
\end{aligned}
\end{equation}
which, for the two-dimensional case, can be written as
    \begin{equation*}
\begin{aligned}
  \bU_{n2}&=\bU_n+\tau\varphi_1(\tau A^\bu_1)
  \left(
  A_1^\bu\bU_n+\bU_n(A_2^\bu)^{\sf T}+
  \bG^\bu(\bU_n,\bV_n)\right)\left(\varphi_1(\tau A_2^\bu)\right)^{\sf T},\\
  \bV_{n2}&=\bV_n+\tau\varphi_1(\tau A^\bv_1)
  \left(  A_1^\bv\bV_n+\bV_n(A_2^\bv)^{\sf T}+
  \bG^\bv(\bU_n,\bV_n)\right)\left(\varphi_1(\tau A_2^\bv)\right)^{\sf T},\\
  \bU_{n+1}&=\bU_{n2}+2\tau \varphi_2(\tau A_1^\bu)(\bG^\bu(\bU_{n2},\bV_{n2})
  -\bG^\bu(\bU_n,\bV_n))\left(\varphi_2(\tau A_2^\bu)\right)^{\sf T},\\
    \bV_{n+1}&=\bV_{n2}+2\tau \varphi_2(\tau A_1^\bv)(\bG^\bv(\bU_{n2},\bV_{n2})
  -\bG^\bv(\bU_n,\bV_n))\left(\varphi_2(\tau A_2^\bv)\right)^{\sf T}.
\end{aligned}
\end{equation*}
{\color{black}
  Notice that the ETD2RKds scheme~\eqref{eq:etd2rk_d_split} is second order
  accurate in time, as can be easily seen by comparing it with the original
  ETD2RK method and using formulas~\eqref{eq:phisplitcomp}.
The expected numerical rate of convergence of the scheme will also be verified
later on in the numerical examples, see the tables in Section~\ref{sec:numexp}.}
{\color{black}Clearly,
  once the relevant small sized matrix functions have been computed,
  the realization of a time step requires four Tucker operators and two
  actions of matrices in Kronecker form. Hence, its computational
  cost mainly depends on the number of degrees of freedom $N$, and not
  on the time step size~$\tau$.}
Finally, remark that
scheme~\eqref{eq:etd2rk_d_split} can be straightforwardly generalized
to systems of semilinear equations with more than two components.
\section{Numerical examples}\label{sec:numexp}
In this section, we show the effectiveness of the proposed
directional split exponential integrator ETD2RKds
by applying it to several two-dimensional and
three-dimensional models of great interest.
All the problems are discretized in space by second order centered
finite differences, {\color{black}both for the diffusion and the advection
terms,} on a grid of equispaced points with
$n_1=\ldots=n_d=n$ (for a total number of degrees of freedom equal to $N=n^d$).
{\color{black}
The usage of centered finite differences for the advection
  terms is justified by the fact that  the cell P\'eclet numbers are much
  less than one in all the relevant examples.
  The discretization of the homogeneous
  Neumann boundary conditions is directly embedded into the matrices.
  Notice that this choice gives rise to tridiagonal
  matrices $A_\mu\in\RR^{n\times n}$.
We remark, however, that this structure is not required by our method.}
{\color{black}Concerning the computation of the small-sized matrix 
$\varphi_\ell$ functions needed by
ETD2RKds, since we will use it in a constant time step size scenario, 
we compute them at the beginning of the time integration by
the \verb+phiquad+
function\footnote{Available at
  \url{https://github.com/caliarim/phisplit}, commit
  \texttt{c67abe3}.}
(introduced in Reference~\cite{CC23}). In short, it approximates the
$\varphi_\ell$ functions by employing
the Gauss--Legendre--Lobatto quadrature rule
applied to formula~\eqref{eq:phiell}, in
combination with the modified scaling and squaring
algorithm~\cite{SW09}. The matrix
exponentials needed in this procedure are computed by the built-in MATLAB
function \verb+expm+, that, for nonsymmetric matrices,
implements a variable
order rational Pad\'e algorithm with scaling 
and squaring~\cite{AMH10}. Notice that,
instead of using \verb+phiquad+,
other techniques could also be employed,
such as direct Pad\'e approximation
of the $\varphi_\ell$ functions (as it is done
in the function \verb+phipade+ of Reference~\cite{BSW05}).}
{\color{black}As already observed in Reference~\cite{CCEOZ22},
also for the examples presented here
  this initial phase has a
  negligible computational cost with respect to the total integration time,
  since the size of the directional matrices $A_\mu$ is small.
  In fact, the maximum size of
  exponential-like matrix functions computed is $200\times 200$
  in the two-dimensional DIB model of Section~\ref{sec:dib}.}

We compare our scheme with other popular exponential or implicit integrators,
developed both for the two-dimensional and the three-dimensional case.
The exponential schemes are employed with constant time step size.
All the methods are implemented in MATLAB language, and
the arising matrix functions or linear systems are computed
by \emph{direct} methods (that is, without using iterative procedures).
In particular, we consider the following schemes.
\begin{itemize}
\item {\color{black}An implementation of the Lawson2b method~\eqref{eq:law2bten}
  in tensor formulation,
  in which the needed small sized matrix exponentials
  are computed, once and for all at the beginning of the time integration,
  by the built-in MATLAB function \verb+expm+. As for the ETD2RKds scheme,
  this initial phase has a negligible computational cost. Notice also
  that this scheme has no directional splitting error.}
\item {\color{black}An implementation of the standard
  ETD2RK method in which the action of the relevant
  matrix functions is realized in tensor formulation by
  the \textsc{phiks}
  routine\footnote{Available at \url{https://github.com/caliarim/phiks},
    commit \texttt{97ae469}.}
  (see end of Section~\ref{sec:expint} and
  Reference~\cite{CCZ23} for a detailed explanation).
  The required small sized matrix exponential
  are computed using the function \verb+expm+.
  The scaling parameter and the number of quadrature nodes are
  automatically chosen by \textsc{phiks}, based on the input
  tolerance. In fact, the tolerances range
  from $5\mathrm{e}{-6}$ to $1\mathrm{e}{-3}$,
  depending on the specific numerical example.}
  \item {\color{black}An implementation of the Exponential Time Differencing
    Real Distinct Poles Integrating
    Factor method (ETD-RDP-IF, see~References~\cite{AAKW20,M22}
    and the accompanying MATLAB
    software\footnote{Available at
      \url{https://github.com/kleefeld80/ETDRDPIF}, commit \texttt{2647b6e}.})
      for the solution of the vector ODE system~\eqref{eq:odes}.
    The exponential functions within the predictor and corrector terms
    are approximated by a first-order Pad\'e expansion
    and a second-order rational approximation with simple real distinct poles.
    After dimensional splitting, thanks to the chosen space discretization,
    the arising linear systems involving the large matrices
$A_{\otimes\mu}\in\RR^{N\times N}$ (that have three nonzero diagonals)
    are solved by an
    adapted Thomas algorithm. We remark that this
    time integrator could benefit
  from a three-core parallel implementation, which we do not consider here.}
\item The MATLAB solver \verb+ode23tb+, which implements
  a variable step size diagonally implicit Runge--Kutta pair 2(3)
  (DIRK23, see Reference~\cite{HS96}) {\color{black}for the
    vector system~\eqref{eq:odes}},
  suggested for stiff problems with ``crude error tolerances'',
  fed with the exact Jacobian of the system.
\end{itemize}
Moreover, we also test the performance of the MATLAB solver \verb+ode23+,
which implements
a variable time step explicit Runge--Kutta pair 3(2)
(RK32, see~Reference~\cite{BS89}) {\color{black}for the
  ODE system~\eqref{eq:odes}.}
This method is considered in the comparisons to be sure that the
examples, with
the selected parameters and space discretizations, are in a stiff regime in
which explicit methods do suffer from a time step size restriction.
{\color{black}As a consequence, we expect that
  the variable step size mechanism will shrink the time step size
  to make the scheme operate in its stability region. Hence,
  the method will be able to reach different accuracies
  with approximately
the same  number of time steps and thus of
computational load (see, for instance, the first experiment of
Section~\ref{sec:Schnakenberg_2d} for more details).}
Whenever calling the built-in MATLAB integrators
{\color{black}we always use sparse
matrices and}, since we are just interested
in the solution at final times, we employ a
proper output function (through the option \verb+OutputFcn+) in
order not to save the solutions at intermediate steps and hence waste
memory. {\color{black}
  We believe that running the experiments  with
  well-established built-in MATLAB integrators is valuable
  and can constitute a
  common and useful benchmark for the community working in the field of
  ADR systems using this
  popular scientific language.}

In all the numerical {\color{black}examples}, we first
compare the performances
of the methods in reaching a common range
of accuracies with respect to a reference
solution computed with ETD2RKds and a sufficiently large number of time steps.
We therefore select, for the different methods,
some sequences of time steps (for the fixed
time step size implementations) or input
tolerances (for the variable time step size
implementations), and measure the overall
needed wall-clock time. The error is computed as the 2-norm
of the relative Frobenius norm of the solutions $\bU$ and $\bV$
at final time. {\color{black}First of all, this experiment is useful to
  check that all the constant time step size methods exhibit the expected
  rate of convergence. In addition, we compute the computational load
  of a single step of our proposed exponential scheme.
Then, in a second experiment, we run ETD2RKds up
to a larger final time, plot the corresponding component $\bU$,
and report the wall-clock time.
Moreover, for the systems of diffusion--reaction equations that lead
to the formation of a stationary Turing pattern, we show the evolution dynamics
by reporting the discretized spatial mean of $u(t_n,\bx)$
\begin{equation*}
  \langle \bU_n\rangle
  \approx  \frac{1}{\lvert\Omega\rvert}\int_{\Omega}u(t_n,\bx)d\bx,
\end{equation*}
and the time increment $\lVert \bU_{n+1}-\bU_n\rVert_\mathrm{F}$ in the
Frobenius norm, as it is typically done in the literature (see, for
instance, References~\cite{AMS23,DASS20}).}

All the experiments are performed on an
Intel\textsuperscript{\textregistered} Core\textsuperscript{\texttrademark} i7-10750H CPU
with six physical cores and
16GB of RAM. As a software, we use MathWorks
MATLAB\textsuperscript{\textregistered} R2022a.
All the codes to reproduce the examples can be
found in a maintained GitHub
repository\footnote{Available at \url{https://github.com/cassinif/ExpADRds}.}.
\subsection{Two-dimensional Schnakenberg model}\label{sec:Schnakenberg_2d}
We consider the Schnakenberg model
(see References~\cite{S79,AMS23,DASS20,BTMSC23})
\begin{equation}\label{eq:Schnakenberg_2d}
  \left\{
  \begin{aligned}
    \partial_t u&=
    \delta^u\Delta u+\rho(a^u-u+u^2v),\\
    \partial_t v&=
    \delta^v \Delta v+\rho (a^v-u^2v),
  \end{aligned}\right.
\end{equation}
in the spatial domain $\Omega=[0,1]^2$. The unknowns $u$ and $v$
represent two chemical concentrations in autocatalytic reactions.
The parameters, taken from Reference~\cite{AMS23}, are
$\delta^u=1$, $\delta^v=10$, $\rho=1000$, $a^u=0.1$, and $a^v=0.9$.
The equilibrium $(u_\rme,v_\rme)=(a^u+a^v,a^v/(a^u+a^v)^2)$
is susceptible
of Turing instability. The initial data are
$u_0=u_\rme+10^{-5}\cdot\mathcal{U}(0,1)$ and
$v_0=v_\rme+10^{-5}\cdot\mathcal{U}(0,1)$,
{\color{black} where here and throughout the experiments
$\mathcal{U}(0,1)$ denotes the uniformly distributed random variable
in $(0,1)$.
The MATLAB random generator seed is set to the value $0$.}
The spatial domain is discretized with
a grid of $N=150^2$ points. The final
simulation time is $T=0.25$.
{\color{black}
The detailed outcome of the experiment is reported in Table~\ref{tab:schnak_2d},
in which we also indicate the numbers of time steps (or input tolerances) and
the obtained numerical order of convergence for the different integrators.
In particular, concerning
the MATLAB variable step size Runge--Kutta integrators,
in this and in all the subsequent examples
the options \verb+AbsTol+ and \verb+RelTol+ are set to the
values of the
reported tolerances.
The results are also graphically presented in a precision diagram 
in Figure~\ref{fig:schnak_2d}.}
\begin{table}[htb!]
  \centering
  {\small
  \begin{tabular}{llllllll}
    \hline
    \multicolumn{4}{l}{ETD2RKds} & \multicolumn{4}{l}{ETD2RK}\\
        steps & time (s) & error & order &     steps & time (s) & error & order\\
    \hline
    3000 & 3.82 & $1.78\rme{-3}$ & ---  & 3000 & 176.98 & $1.65\rme{-3}$ & ---     \\
    4000 & 5.17 & $1.01\rme{-3}$ & 1.98 & 4000 & 227.38 & $9.33\rme{-4}$ & 1.97 \\
    5000 & 6.43 & $6.42\rme{-4}$ & 2.02 & 5000 & 282.61 & $5.95\rme{-4}$ & 2.02 \\
    6000 & 7.73 & $4.41\rme{-4}$ & 2.06 & 6000 & 334.67 & $4.09\rme{-4}$ & 2.06 \\
    \hline
    \hline
    \multicolumn{4}{l}{Lawson2b} & \multicolumn{4}{l}{ETD-RDP-IF}\\
        steps & time (s) & error & order &     steps & time (s)& error & order\\
    \hline
    14000 & 10.56 & $3.58\rme{-3}$ & ---  & 14000 & 97.40 & $1.99\rme{-3}$ & --- \\
    18000 & 13.33 & $2.16\rme{-3}$ & 2.01 & 18000 & 125.23 & $1.21\rme{-3}$ & 2.00\\
    22000 & 16.07 & $1.44\rme{-3}$ & 2.02 & 22000 & 152.67 & $8.05\rme{-4}$ & 2.02\\
    26000 & 19.01 & $1.03\rme{-3}$ & 2.03 & 26000 & 180.81 & $5.71\rme{-4}$ & 2.05\\
    \hline
    \hline
    \multicolumn{4}{l}{DIRK23} & \multicolumn{4}{l}{RK32}\\
        tolerance & time (s) & error & &     tolerance & time (s) & error & \\
    \hline
    $1\mathrm{e}{-6}$ & 10.58 & $5.12\rme{-3}$ & & $1\mathrm{e}{-2}$ & 366.16 & $4.37\rme{-3}$ & \\
    $5\mathrm{e}{-7}$ & 11.57 & $2.91\rme{-3}$ & & $8\mathrm{e}{-3}$ & 371.76 & $2.93\rme{-3}$ & \\
    $1\mathrm{e}{-7}$ & 16.40 & $9.37\rme{-4}$ & & $4\mathrm{e}{-3}$ & 373.26 & $1.10\rme{-3}$ & \\
    $5\mathrm{e}{-8}$ & 18.66 & $5.54\rme{-4}$ & & $1\mathrm{e}{-3}$ & 377.94 & $6.47\rme{-4}$ & \\
    \hline
  \end{tabular}}%
  \caption{{\color{black}Number of time steps (or input tolerance), wall-clock time (in seconds),
  relative error at final time, and observed numerical order of convergence
      for the solution of the 2D Schnakenberg model~\eqref{eq:Schnakenberg_2d} 
      up to $T=0.25$ with different integrators.
      See also Figure~\ref{fig:schnak_2d} for a graphical representation.}}
  \label{tab:schnak_2d}
\end{table}
\begin{subfigures}
  \begin{figure}[htb!]
  \centering
%
%
\definecolor{mycolor1}{rgb}{1.00000,0.00000,1.00000}%
\definecolor{mycolor2}{rgb}{0.00000,1.00000,1.00000}%
\begin{tikzpicture}

\begin{axis}[%
width=2.5in,
height=2in,
at={(0.769in,0.477in)},
scale only axis,
xmode=log,
xmin=2,
xmax=1000,
xminorticks=true,
ymode=log,
ymin=0.00008,
ymax=0.01,
xlabel = {Wall-clock time (s)},
ylabel = {Relative error},
yminorticks=true,
axis background/.style={fill=white},
legend style={at={(0.05,0.06)}, anchor=south west, legend cell align=left, align=left, draw=white!15!black, font=\scriptsize},
legend columns = 2
]
\addplot [color=mycolor1,line width=1.5pt, mark size = 3pt, mark=+, mark options={solid, mycolor1}]
  table[row sep=crcr]{%
3.818152	0.00178203082851037\\
5.171972	0.00100801054590846\\
6.429186	0.000642364808667461\\
7.72870600000001	0.000441072294356976\\
};
\addlegendentry{ETD2RKds}

\addplot [color=blue, line width=1.5pt, mark size = 3pt,mark=o, mark options={solid, blue}]
  table[row sep=crcr]{%
176.975237	0.00164538403126529\\
227.384505	0.000933216679207764\\
282.607541	0.000595260399777974\\
334.669068	0.000408706249872652\\
};
\addlegendentry{ETD2RK}

\addplot [color=red,line width=1.5pt, mark size = 3pt, mark=triangle, mark options={solid, rotate=270, red}]
  table[row sep=crcr]{%
10.55984	0.00358441428281133\\
13.327032	0.0021643311007106\\
16.074708	0.00144378482262667\\
19.009404	0.00102839256819669\\
};
\addlegendentry{Lawson2b}


\addplot [color=green, line width=1.5pt, mark size = 3pt, mark=asterisk, mark options={solid, green}]
  table[row sep=crcr]{%
97.40	0.00199458392217813\\
125.23	0.00120698622678295\\
152.67	0.000804511002692113\\
180.81	0.000571389367311821\\
};
\addlegendentry{ETD-RDP-IF}

\addplot [color=mycolor2,line width=1.5pt, mark size = 3pt, mark=triangle, mark options={solid, mycolor2}]
  table[row sep=crcr]{%
10.575728	0.00512055647843696\\
11.565102	0.00291275422053487\\
16.395358	0.000936721234779732\\
18.656363	0.000553522235380911\\
};
\addlegendentry{DIRK23}

\addplot [color=black,line width=1.5pt, mark size = 3pt, mark=triangle, mark options={solid, rotate=90, black}]
  table[row sep=crcr]{%
366.164689	0.004365215747080\\
371.76735	 0.002926081526388\\
373.26104	 0.001099095908898\\
377.93894	 6.465162417369081e-04\\
};
\addlegendentry{RK32}

\end{axis}

\end{tikzpicture}%
  \caption{Results for the simulation of the 2D Schnakenberg
    model~\eqref{eq:Schnakenberg_2d} with
  $N=150^2$ spatial discretization points. The number
  of time steps (or input tolerance) for each integrator is reported in
  Table~\ref{tab:schnak_2d}. The final simulation time is $T=0.25$.}
  \label{fig:schnak_2d}
\end{figure}
\begin{figure}[htb!]
  \centering
  \input{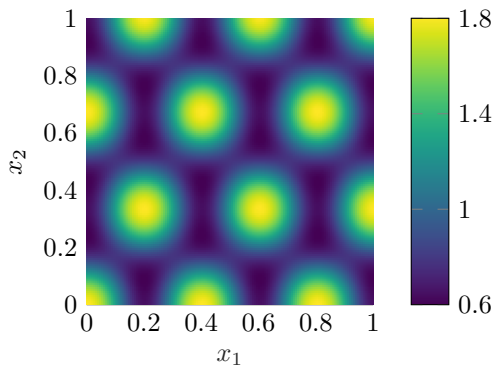}
  \caption{Turing pattern ($u$ component)
    for 2D Schnakenberg model~\eqref{eq:Schnakenberg_2d}
    obtained at final time
    $T=2$ with $N=150^2$ space discretization points with ETD2RKds.
The time step size employed is {\color{black}$\tau=5\mathrm{e}{-4}$}.
The simulation wall-clock time is {\color{black}5} seconds.}
  \label{fig:schnak2d_patt}
\end{figure}
\begin{figure}
  \centering
  \input{img/schnak2D_mean.tex}
  \input{img/schnak2D_increment.tex}
  \caption{{\color{black} Indicators for the time dynamics of 
      2D Schnakenberg model~\eqref{eq:Schnakenberg_2d} solved up to final time
      $T=2$ with ETD2RKds
    and time step size $\tau=5\mathrm{e}{-4}$. The top plot refers to the spatial mean
    $\langle \bU_n \rangle$ while the bottom plot depicts the time increment
    $\lVert\bU_{n+1}-\bU_n \rVert_\mathrm{F}$}.}
  \label{fig:Schnak2D_indicators}
\end{figure}
\end{subfigures}%

{\color{black}First of all, we observe that all the
constant time step size methods exhibit the expected rate
of convergence (see the relevant columns in Table~\ref{tab:schnak_2d}).}
Then, concerning Figure~\ref{fig:schnak_2d},
the plot is clearly divided into two parts. On the right we have the
more expensive methods.
{\color{black}As expected, the RK32 method performs poorly in terms
  of computational time. Moreover, we observe that it requires
  about  the same
computational load for the four different input tolerances.
In fact, the numbers of time steps performed by the RK32 method
are 176797, 176798, 176800, and 176802, respectively.
This behavior is expected, since the scheme is explicit and
is forced to take a large
number of time steps due to stability constraints.
Nevertheless, due to the automatic time step selection, slightly distinct
sequences of time steps are produced, and hence different
final accuracies are obtained.}
{\color{black}
The ETD-RDP-IF scheme (specifically designed for these type of
models) performs slightly better than the 
ETD2RK method, both schemes
being however not competitive with the remaining methods.
On the left hand side of the plot,
we find the implicit Runge--Kutta method,
which is more than ten times faster
than the explicit one, and with performances similar to the Lawson2b method.
Finally,
an overall neat advantage is obtained by applying the proposed method ETD2RKds,
which requires an average wall-clock time of $1.3\mathrm{e}{-3}$ seconds per time step.}

In the second experiment,
we perform a simulation with ETD2RKds up to final time
$T=2$ with {\color{black}4000} time steps.
The observed
pattern is in agreement with
that already reported in the literature (compare Figure~\ref{fig:schnak2d_patt}
and Reference~\cite[Fig.~8(a)]{AMS23}). {\color{black}
 The overall wall-clock time of the simulation
is roughly 5 seconds,
  which corresponds
  to the just mentioned average time step cost multiplied by the
  number of steps, as expected.}
  {\color{black}Finally, we report in Figure~\ref{fig:Schnak2D_indicators} two useful indicators
    for the dynamics of the system, i.e., the spatial mean and the time
    increments.
  In particular, it is clearly visible the distinction between the initial 
  reactivity phase (up to about $t=0.25$) and the stabilization to a 
  stationary pattern.}
\subsection{Two-dimensional FitzHugh--Nagumo model}\label{sec:FitzHughNagumo}
We consider the FitzHugh--Nagumo model (see Reference~\cite{AMS23})
\begin{equation}\label{eq:FitzHughNagumo}
  \left\{
  \begin{aligned}
    \partial_t u&=
    \delta^u\Delta u+\rho(-u(u^2-1)-v),\\
    \partial_t v&=
    \delta^v \Delta v+\rho a^v_1(u-a_2^vv),
  \end{aligned}\right.
\end{equation}
in $\Omega=[0,\pi]^2$. The model describes the flow of an electric current
through a nerve fiber. The unknowns $u$ and $v$ represent
the electric potential and the recovery variable, respectively.
The parameters, taken from Reference~\cite{AMS23},
are $\delta^u=1$, $\delta^v=42.1887$,  $\rho=65.731$, $a_1^v=11$, and
$a_2^v=0.1$.
With these choices, the equilibrium $(u_\rme,v_\rme)=(0,0)$ is susceptible
of Turing instability.
The initial solutions are $u_0=10^{-3}\cdot\mathcal{U}(0,1)$
and $v_0=10^{-3}\cdot\mathcal{U}(0,1)$
{\color{black}with MATLAB random generator seed set to the value $0$.}
The spatial domain is discretized by a grid of $N=100^2$ points.
In the first experiment we simulate
up to $T=10$.
The numbers of time steps (or input tolerances) for each integrator, as well
as the detailed outcome of the experiment, can be found in
Table~\ref{tab:fhn}. The results are graphically presented in Figure~\ref{fig:fhn_2d}. 
\begin{table}[htb!]
  \centering
  {\small
  \begin{tabular}{llllllll}
    \hline
    \multicolumn{4}{l}{ETD2RKds} & \multicolumn{4}{l}{ETD2RK}\\
        steps & time (s) & error & order &     steps & time (s) & error & order\\
    \hline
    20000 & 13.50 & $9.16\rme{-3}$ & ---  & 20000 & 499.28 & $8.34\rme{-3}$ & ---     \\
    22500 & 17.66 & $7.36\rme{-3}$ & 1.85 & 22500 & 560.11 & $6.73\rme{-3}$ & 1.83 \\
    25000 & 19.35 & $6.04\rme{-3}$ & 1.88 & 25000 & 627.55 & $5.53\rme{-3}$ & 1.86 \\
    27500 & 21.21 & $5.03\rme{-3}$ & 1.91 & 27500 & 686.38 & $4.61\rme{-3}$ & 1.89 \\
    \hline
    \hline
    \multicolumn{4}{l}{Lawson2b} & \multicolumn{4}{l}{ETD-RDP-IF}\\
        steps & time (s) & error & order &     steps & time (s)& error & order\\
    \hline
    600000 & 277.90 & $3.86\rme{-2}$ & ---  & 600000 & 1794.51 & $2.57\rme{-2}$ & --- \\
    675000 & 307.72 & $3.04\rme{-2}$ & 2.02 & 675000 & 2025.12 & $2.03\rme{-2}$ & 2.01\\
    750000 & 337.09 & $2.46\rme{-2}$ & 2.02 & 750000 & 2263.78 & $1.64\rme{-2}$ & 2.02\\
    825000 & 362.99 & $2.03\rme{-2}$ & 2.02 & 825000 & 2476.61 & $1.35\rme{-2}$ & 2.02\\
    \hline
    \hline
    \multicolumn{4}{l}{DIRK23} & \multicolumn{4}{l}{RK32}\\
        tolerance & time (s) & error & &     tolerance & time (s) & error & \\
    \hline
    $5\mathrm{e}{-7}$ & 33.59 & $1.33\rme{-2}$ & & $8\mathrm{e}{-7}$ & 1548.62 & $2.29\rme{-2}$ & \\
    $1\mathrm{e}{-7}$ & 47.48 & $7.89\rme{-3}$ & & $6\mathrm{e}{-7}$ & 1601.73 & $1.48\rme{-2}$ & \\
    $8\mathrm{e}{-8}$ & 50.95 & $6.39\rme{-3}$ & & $4\mathrm{e}{-7}$ & 1641.75 & $6.16\rme{-3}$ & \\
    $6\mathrm{e}{-8}$ & 53.07 & $5.70\rme{-3}$ & & $2\mathrm{e}{-7}$ & 1654.41 & $2.53\rme{-3}$ & \\
    \hline
  \end{tabular}}%
  \caption{{\color{black}Number of time steps (or input tolerance), wall-clock time (in seconds),
  relative error at final time, and observed numerical order of convergence
      for the solution of the 2D FitzHugh--Nagumo model~\eqref{eq:FitzHughNagumo} 
      up to $T=10$ with different integrators.
      See also Figure~\ref{fig:fhn_2d} for a graphical representation.}}
  \label{tab:fhn}
\end{table}
\begin{subfigures}
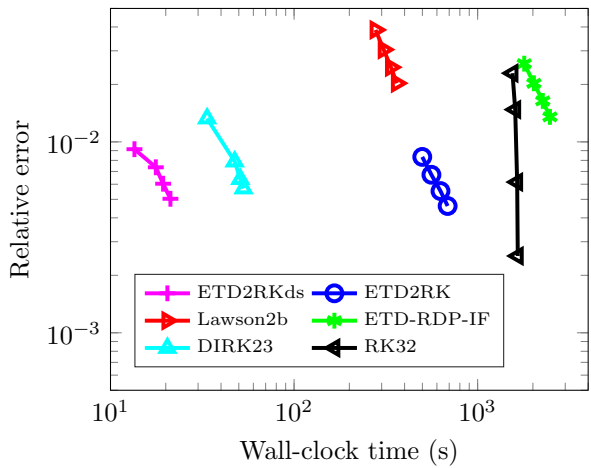
\begin{figure}[htb!]
  \centering
%
%
\definecolor{mycolor1}{rgb}{1.00000,0.00000,1.00000}%
\definecolor{mycolor2}{rgb}{0.00000,1.00000,1.00000}%
\begin{tikzpicture}

\begin{axis}[%
width=2.5in,
height=2in,
at={(0.769in,0.477in)},
scale only axis,
xmode=log,
xmin=10,
xmax=4000,
xminorticks=true,
ymode=log,
ymin=0.0005,
ymax=0.05,
xlabel = {Wall-clock time (s)},
ylabel = {Relative error},
yminorticks=true,
axis background/.style={fill=white},
legend style={at={(0.05,0.06)}, anchor=south west, legend cell align=left, align=left, draw=white!15!black, font=\scriptsize},
legend columns = 2
]
\addplot [color=mycolor1,line width=1.5pt, mark size = 3pt, mark=+, mark options={solid, mycolor1}]
  table[row sep=crcr]{%
13.504601	0.00915542026217866\\
17.656152	0.0073613901297128\\
19.353818	0.00603818597888116\\
21.205187	0.00503429722077697\\
};
\addlegendentry{ETD2RKds}

\addplot [color=blue, line width=1.5pt, mark size = 3pt,mark=o, mark options={solid, blue}]
  table[row sep=crcr]{%
499.283052	0.00834257646264861\\
560.110838	0.00672523569620684\\
627.547039	0.00552707248369243\\
686.383449	0.00461485879407511\\
};
\addlegendentry{ETD2RK}

\addplot [color=red,line width=1.5pt, mark size = 3pt, mark=triangle, mark options={solid, rotate=270, red}]
  table[row sep=crcr]{%
277.902688	0.038604904035945\\
307.723031	0.0304319653858592\\
337.094255	0.0245992481705008\\
362.992432	0.0202899769415397\\
};
\addlegendentry{Lawson2b}


\addplot [color=green, line width=1.5pt, mark size = 3pt, mark=asterisk, mark options={solid, green}]
  table[row sep=crcr]{%
1794.51	0.025690110857349\\
2025.12	0.020270247132495\\
2263.78	0.0163925611503285\\
2476.61	0.0135222222496991\\
};
\addlegendentry{ETD-RDP-IF}

\addplot [color=mycolor2,line width=1.5pt, mark size = 3pt, mark=triangle, mark options={solid, mycolor2}]
  table[row sep=crcr]{%
33.589773	0.013265638263018\\
47.480417	0.007893035952687\\
50.951479	0.006386653705674\\
53.065652	0.005703427262398\\
};
\addlegendentry{DIRK23}

\addplot [color=black,line width=1.5pt, mark size = 3pt, mark=triangle, mark options={solid, rotate=90, black}]
  table[row sep=crcr]{%
1.5486e+03	2.2894e-02\\
1.6017e+03	1.4762e-02\\
1.6417e+03	6.1585e-03\\
1.6544e+03	2.5258e-03\\
};
\addlegendentry{RK32}

\end{axis}

\end{tikzpicture}%
  \caption{Results for the simulation of 2D FitzHugh--Nagumo
  model~\eqref{eq:FitzHughNagumo} with
  $N=100^2$ spatial discretization points. The number
  of time steps (or input tolerance) for each integrator is reported in
  Table~\ref{tab:fhn}. The final simulation time is $T=10$.}
  \label{fig:fhn_2d}
\end{figure}
\begin{figure}[htb!]
  \centering
  \input{img/Upatt_fhn_2d_100.tex}
  \caption{Turing pattern ($u$ component) for 2D
    FitzHugh--Nagumo model~\eqref{eq:FitzHughNagumo}
    obtained at final time
    $T=50$ with $N=100^2$ space discretization points with ETD2RKds. The
    time step size employed is {\color{black}$\tau=1.67\mathrm{e}{-3}$}. The simulation
    wall-clock time is {\color{black}22} seconds.}
  \label{fig:fhn_patt}
\end{figure}
\begin{figure}
  \centering
  \input{img/fhn2D_mean.tex}
  \input{img/fhn2D_increment.tex}
  \caption{{\color{black} Indicators for the time dynamics of 
      2D FitzHugh--Nagumo model~\eqref{eq:FitzHughNagumo} solved up to final
      time $T=50$ with ETD2RKds and
    time step size $\tau=1.67\mathrm{e}{-3}$. The top plot refers to the spatial mean
    $\langle \bU_n \rangle$ while the bottom plot depicts the time increment
    $\lVert\bU_{n+1}-\bU_n \rVert_\mathrm{F}$}.}
  \label{fig:Fhn2D_indicators}
\end{figure}
\end{subfigures}

First of
all, we notice that the Lawson2b and ETD-RDP-IF methods require many more
time steps than the exponential Runge--Kutta methods, and still they
cannot
reach the same accuracies. Their wall-clock time is in fact at least
one order of magnitude larger than the most efficient methods.
On the other hand, as in the previous example,
the ETD2RKds method is not heavily affected by
the directional splitting error. Indeed, it reaches almost the same
accuracies obtained by the ETD2RK method but with a much smaller computational
time. {\color{black}In this case, the average wall-clock time per time step of ETD2RKds
is $7.5\mathrm{e}{-4}$ seconds.}

In the second experiment, we set the final simulation time to $T=50$
and integrate
the problem with the ETD2RKds method with {\color{black}30000} time steps.
The overall computational time is about {\color{black}22} seconds.
Again, the obtained pattern agrees with that reported in the literature
(compare Figure~\ref{fig:fhn_patt} and Reference~\cite[Fig.~4(a)]{AMS23}),
i.e., a square pattern corresponding to the cosine modes $(4,4)$.
{\color{black}Also, as expected both the indicators reported in
Figure~\ref{fig:Fhn2D_indicators} tend to zero.}
\subsection{Three-dimensional
  FitzHugh--Nagumo model}\label{sec:FitzHughNagumo_3d}
We consider again the FitzHugh--Nagumo model
\begin{equation}\label{eq:FitzHughNagumo_3d}
  \left\{
  \begin{aligned}
    \partial_t u&=
    \delta^u\Delta u+\rho(-u(u^2-1)-v),\\
    \partial_t v&=
    \delta^v \Delta v+\rho a_1^v(u-a_2^vv),
  \end{aligned}\right.
\end{equation}
but now in a three-dimensional domain
$\Omega=[0,\pi]^3$. 
The parameters
$\delta^u=1$, $\delta^v=42.1887$,  $\rho=24.649$, $a_1^v=11$, and $a_2^v=0.1$
were obtained by following the analysis performed
in Reference~\cite{GLRS19} to achieve a stationary square pattern
with modes $(2,2,2)$.
With this choice of parameters,
the equilibrium $(u_\rme,v_\rme)=(0,0)$ is indeed susceptible
of Turing instability.
The initial solutions are $u_0=10^{-3}\cdot\mathcal{U}(0,1)$
and $v_0=10^{-3}\cdot\mathcal{U}(0,1)$,
{\color{black}with MATLAB random generator seed set to the value $0$.}
The spatial domain is discretized by a grid of $N=64^3$ points.
We simulate up to the final time $T=10$
with the number of time steps given in
Table~\ref{tab:FitzHughNagumo_3d}. The results, reported in the table, 
are also graphically depicted in
Figure~\ref{fig:FitzHughNagumo_3d}.
\begin{table}[htb!]
  \centering
  {\small
  \begin{tabular}{llllllll}
    \hline
    \multicolumn{4}{l}{ETD2RKds} & \multicolumn{4}{l}{ETD2RK}\\
        steps & time (s) & error & order &     steps & time (s) & error & order\\
    \hline
    12000 & 101.14 & $3.01\rme{-3}$ & ---  & 12000 & 626.16 & $2.67\rme{-3}$ & ---     \\
    14000 & 119.30 & $2.20\rme{-3}$ & 2.03 & 14000 & 718.69 & $1.96\rme{-3}$ & 2.00 \\
    16000 & 138.89 & $1.67\rme{-3}$ & 2.05 & 16000 & 816.36 & $1.49\rme{-3}$ & 2.03 \\
    18000 & 154.95 & $1.31\rme{-3}$ & 2.07 & 18000 & 889.12 & $1.17\rme{-3}$ & 2.06 \\
    \hline
  \end{tabular}
  \begin{tabular}{llll}
    \hline
    \multicolumn{3}{l}{Lawson2b}\\
        steps & time (s) & error & order \\
    \hline
    250000 & 1429.71 & $1.13\rme{-2}$ & --- \\
    300000 & 1729.65 & $7.86\rme{-3}$ & 2.00\\
    350000 & 2013.62 & $5.77\rme{-3}$ & 2.00\\
    400000 & 2320.22 & $4.42\rme{-3}$ & 2.01\\
    \hline
  \end{tabular}}%
  \caption{{\color{black}Number of time steps, wall-clock time (in seconds),
  relative error at final time, and observed numerical order of convergence
      for the solution of the 3D FitzHugh--Nagumo model~\eqref{eq:FitzHughNagumo_3d} 
      up to $T=10$ with different integrators.
      The DIRK23 method interrupted due to
      excessive memory requirements with tolerance $1\mathrm{e}{-1}$.
      The RK32 and the ETD-RDP-IF methods did not output a solution
      within $10^4$ seconds.
      See also Figure~\ref{fig:FitzHughNagumo_3d} for a graphical representation.}}
  \label{tab:FitzHughNagumo_3d}
\end{table}
\begin{subfigures}
\begin{figure}[htb!]
  \centering
%
%
\definecolor{mycolor1}{rgb}{1.00000,0.00000,1.00000}%
\begin{tikzpicture}

\begin{axis}[%
width=2.5in,
height=2in,
at={(0.769in,0.477in)},
scale only axis,
xmode=log,
xmin=40,
xmax=10000,
xminorticks=true,
ymode=log,
ymin=0.001,
ymax=0.02,
xlabel = {Wall-clock time (s)},
ylabel = {Relative error},
yminorticks=true,
axis background/.style={fill=white},
legend style={at={(0.1,0.7)}, anchor=south west, legend cell align=left, align=left, draw=white!15!black, font=\scriptsize},
]
\addplot [color=mycolor1,line width=1.5pt, mark size = 3pt, mark=+, mark options={solid, mycolor1}]
  table[row sep=crcr]{%
101.142809	0.00301007064517505\\
119.296966	0.00220144307325069\\
138.889036	0.00167373583396444\\
154.948518	0.00131096712188644\\
};
\addlegendentry{ETD2RKds}

\addplot [color=blue, line width=1.5pt, mark size = 3pt,mark=o, mark options={solid, blue}]
  table[row sep=crcr]{%
626.163961	0.00266506822397421\\
718.690042	0.00195927768365453\\
816.35759	0.00149409782391268\\
889.122835	0.00117195916241373\\
};
\addlegendentry{ETD2RK}

\addplot [color=red,line width=1.5pt, mark size = 3pt, mark=triangle, mark options={solid, rotate=270, red}]
  table[row sep=crcr]{%
1.4297e3	1.1321e-2\\
1.7296e3	7.8615e-3\\
2.0136e3	5.7731e-3\\
2.3202e3	4.4161e-3\\
};
\addlegendentry{Lawson2b}



\end{axis}

\end{tikzpicture}%
  \caption{Results for the simulation of the 3D FitzHugh--Nagumo
  model~\eqref{eq:FitzHughNagumo_3d} with
  $N=64^3$ spatial discretization points. The number
  of time steps for each integrator is reported in
  Table~\ref{tab:FitzHughNagumo_3d}. The final simulation time
  is $T=10$.}
  \label{fig:FitzHughNagumo_3d}
\end{figure}
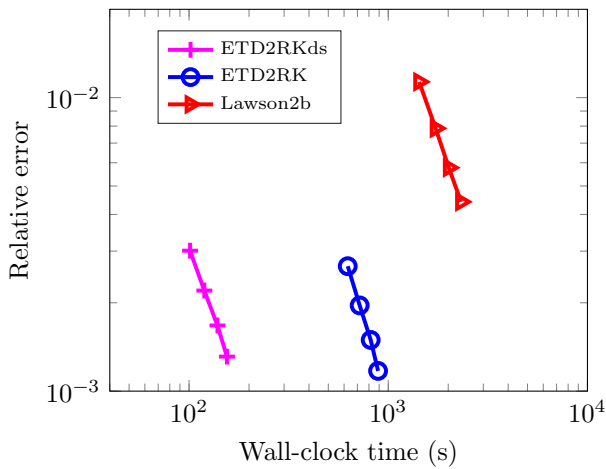
\begin{figure}[htb!]
  \centering
  \input{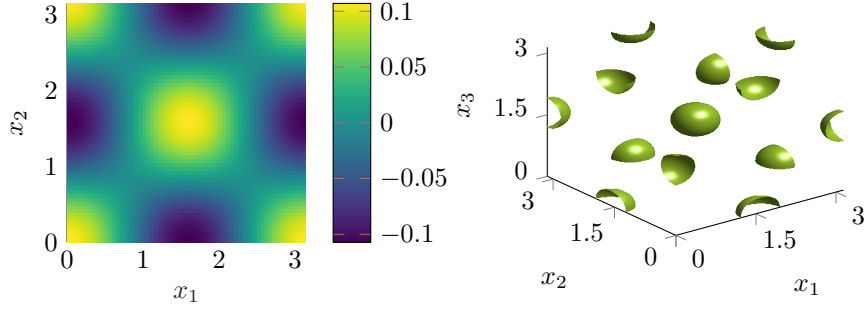}
  \caption{Turing pattern ($u$ component) for 3D FitzHugh--Nagumo
    model~\eqref{eq:FitzHughNagumo_3d}
    obtained at final time
    {\color{black}$T=150$} with $N=64^3$ space discretization points and ETD2RKds. The
    time step size employed is {\color{black}$\tau=6\mathrm{e}{-3}$}.
    The wall-clock simulation time is {\color{black}230} seconds.
    The reported slice (left plot) correspond to $x_3=1.55$
    and the isosurface value (right plot) is $0.08$.}
  \label{fig:FitzHughNagumo_3d_patt}
\end{figure}
\begin{figure}
  \centering
  \input{img/fhn3D_mean.tex}
  \input{img/fhn3D_increment.tex}
  \caption{{\color{black} Indicators for the time dynamics of 
      3D FitzHugh--Nagumo model~\eqref{eq:FitzHughNagumo_3d} solved up to
      final time $T=150$ 
    with  ETD2RKds
    and time step size $\tau=6\mathrm{e}{-3}$. The top plot refers to the spatial mean
    $\langle \bU_n \rangle$ while the bottom plot depicts the time increment
    $\lVert\bU_{n+1}-\bU_n \rVert_\mathrm{F}$}.}
  \label{fig:Fhn3D_indicators}
\end{figure}
\end{subfigures}

First of all, despite  the aforementioned cares in the usage of the
DIRK23 method from the MATLAB ODE suite,
the \verb+ode23tb+ function stopped to run after some seconds
due to too large memory requirements.
Moreover, the RK32 and ETD-RDP-IF methods
did not reach accuracies comparable with the other schemes within $10^4$ seconds.
Among the remaining methods,
ETD2RKds is again the best one, with a speedup of almost one order
of magnitude with respect to the second one, that is the
exponential Runge--Kutta method with the $\varphi_\ell$ functions approximated
by \textsc{phiks}.
{\color{black}The average wall-clock time per time step for ETD2RKds is 
$8.6\mathrm{e}{-3}$ seconds.}

We then simulate the system
up to {\color{black}$T=150$ with 25000} time steps and the ETD2RKds method.  We report in
Figure~\ref{fig:FitzHughNagumo_3d_patt} the obtained pattern {\color{black}and in
Figure~\ref{fig:Fhn3D_indicators} the relevant indicators,} which
agree with the theoretical expectations. The overall wall-clock time
is about {\color{black}230} seconds.

\subsection{Two-dimensional morpho-chemical DIB model}\label{sec:dib}
We consider the morpho-chemical DIB model (see References~\cite{DASS20,AMS23})
\begin{equation}\label{eq:dib}
  \left\{
  \begin{aligned}
    \partial_t u&=
    \delta^u\Delta u+\rho(a^u_1(1-v)u-a^u_2u^3-a^u_3(v-a^u_4)),\\
    \partial_t v&=
    \delta^v \Delta v+\rho(a^v_1(1+a^v_2u)(1-v)(1-a^v_3(1-v))-a^v_4v(1+a^v_5u)(1+a^v_3 v)),
  \end{aligned}\right.
\end{equation}
in $\Omega=[0,20]^2$. The model describes the electrodeposition process for
metal growth. The unknowns $u$ and $v$ represent the morphology of the
metal deposit and its surface chemical composition, respectively.
The parameters, taken from Reference~\cite{AMS23}, are $\delta^u=1$,
$\rho=25/4$,
$a^u_1=10$, $a^u_2=1$, $a_3^u=66$, $a_4^u=0.5$,
$\delta^v=20$, $a_1^v=3$, $a_2^v=2.5$, $a_3^v=0.2$,
\begin{equation*}
  a_4^v=\frac{a_1^v(1-a_4^u)(1-a_3^v+a_3^va_4^u)}{a_4^u(1+a_3^va_4^u)},
\end{equation*}
and $a_5^v=1.5$.
The particular choice of the parameter $a_4^v$ makes the equilibrium
$(u_\rme,v_\rme)=(0,a_4^u)$ susceptible of Turing instability.
The initial conditions are
$u_0=u_\rme+10^{-5}\cdot\mathcal{U}(0,1)$ and
$v_0=v_\rme+10^{-5}\cdot\mathcal{U}(0,1)$.
{\color{black}The MATLAB random generator seed is fixed to the
  value 123.}
The spatial domain is discretized with a grid of $N=200^2$ points.
We first integrate the system up to the final time $T=2.5$ with
a number of time steps, or input tolerance,
as reported in Table~\ref{tab:DIB_ts}.
The detailed results are reported in the table and also presented in the
precision
diagram in Figure~\ref{fig:DIB_cpudiag}.%
\begin{table}[htb!]
  \centering
  {\small
  \begin{tabular}{llllllll}
    \hline
    \multicolumn{4}{l}{ETD2RKds} & \multicolumn{4}{l}{ETD2RK}\\
        steps & time (s) & error & order &     steps & time (s) & error & order\\
    \hline
    1250 & 5.51 & $1.13\rme{-2}$ & ---  & 1250 & 133.63 & $9.48\rme{-3}$ & ---     \\
    1500 & 6.72 & $7.80\rme{-3}$ & 2.04 & 1500 & 146.90 & $6.59\rme{-3}$ & 1.99 \\
    1750 & 7.82 & $5.69\rme{-3}$ & 2.04 & 1750 & 167.88 & $4.84\rme{-3}$ & 2.00 \\
    2000 & 9.13 & $4.34\rme{-3}$ & 2.04 & 2000 & 190.90 & $3.70\rme{-3}$ & 2.01 \\
    \hline
    \hline
    \multicolumn{4}{l}{Lawson2b} & \multicolumn{4}{l}{ETD-RDP-IF}\\
        steps & time (s) & error & order &     steps & time (s)& error & order\\
    \hline
    3750 & 11.52 & $1.62\rme{-2}$ & ---  & 5000 & 62.78 & $2.29\rme{-2}$ & --- \\
    4500 & 13.23 & $1.16\rme{-2}$ & 1.83 & 6000 & 75.46 & $1.62\rme{-2}$ & 1.91\\
    5250 & 16.40 & $8.72\rme{-3}$ & 1.86 & 7000 & 88.46 & $1.20\rme{-2}$ & 1.92\\
    6000 & 18.74 & $6.78\rme{-3}$ & 1.88 & 8000 & 101.04 & $9.29\rme{-3}$ & 1.94\\
    \hline
    \hline
    \multicolumn{4}{l}{DIRK23} & \multicolumn{4}{l}{RK32}\\
        tolerance & time (s) & error & &     tolerance & time (s) & error & \\
    \hline
    $5\mathrm{e}{-5}$ & 15.21 & $2.90\rme{-2}$ & & $9\mathrm{e}{-4}$ & 58.74 & $1.83\rme{-2}$ & \\
    $3\mathrm{e}{-5}$ & 17.04 & $2.04\rme{-2}$ & & $8\mathrm{e}{-4}$ & 59.27 & $1.40\rme{-2}$ & \\
    $2\mathrm{e}{-5}$ & 18.56 & $1.19\rme{-2}$ & & $7\mathrm{e}{-4}$ & 58.62 & $1.07\rme{-2}$ & \\
    $1\mathrm{e}{-5}$ & 23.36 & $5.64\rme{-3}$ & & $6\mathrm{e}{-4}$ & 58.32 & $7.56\rme{-3}$ & \\
    \hline
  \end{tabular}}%
  \caption{{\color{black}Number of time steps (or input tolerance), wall-clock time (in seconds),
  relative error at final time, and observed numerical order of convergence
      for the solution of the 2D DIB model~\eqref{eq:dib} 
      up to $T=2.5$ with different integrators.
      See also Figure~\ref{fig:DIB_cpudiag} for a graphical representation.}}
  \label{tab:DIB_ts}
\end{table}%
\begin{subfigures}
\begin{figure}[htb!]
  \centering
%
%
\definecolor{mycolor1}{rgb}{1.00000,0.00000,1.00000}%
\definecolor{mycolor2}{rgb}{0.00000,1.00000,1.00000}%
\begin{tikzpicture}

\begin{axis}[%
width=2.5in,
height=2in,
scale only axis,
xmode=log,
xmin=4,
xmax=300,
xminorticks=true,
ymode=log,
ymin=3e-3,
ymax=1e-1,
xlabel = {Wall-clock time (s)},
ylabel = {Relative error},
yminorticks=true,
axis background/.style={fill=white},
legend style={at={(0.1,0.7)}, anchor=south west, legend cell align=left, align=left, draw=white!15!black, font=\scriptsize},
legend columns = 2
]
\addplot [color=mycolor1,line width=1.5pt, mark size = 3pt, mark=+, mark options={solid, mycolor1}]
  table[row sep=crcr]{%
5.51	1.130e-02\\
6.72	7.795e-03\\
7.82	5.694e-03\\
9.13	4.337e-03\\
};
\addlegendentry{ETD2RKds}

\addplot [color=blue, line width=1.5pt, mark size = 3pt,mark=o, mark options={solid, blue}]
  table[row sep=crcr]{%
133.63	9.481e-03\\
146.90	6.591e-03\\
167.88	4.839e-03\\
190.90	3.699e-03\\
};
\addlegendentry{ETD2RK}

\addplot [color=red,line width=1.5pt, mark size = 3pt, mark=triangle, mark options={solid, rotate=270, red}]
  table[row sep=crcr]{%
11.52	1.618e-02\\
13.23	1.160e-02\\
16.40	8.715e-03\\
18.74	6.784e-03\\
};
\addlegendentry{Lawson2b}

\addplot [color=green, line width=1.5pt, mark size = 3pt, mark=asterisk, mark options={solid, green}]
  table[row sep=crcr]{%
62.78	2.291e-02\\
75.46	1.618e-02\\
88.46	1.203e-02\\
101.04	9.288e-03\\
};
\addlegendentry{ETD-RDP-IF}

\addplot [color=mycolor2,line width=1.5pt, mark size = 3pt, mark=triangle, mark options={solid, mycolor2}]
  table[row sep=crcr]{%
15.21	2.901e-02\\
17.04	2.037e-02\\
18.56	1.190e-02\\
23.36	5.636e-03\\
};
\addlegendentry{DIRK23}

\addplot [color=black,line width=1.5pt, mark size = 3pt, mark=triangle, mark options={solid, rotate=90, black}]
  table[row sep=crcr]{%
58.74	1.832e-02\\
59.27	1.404e-02\\
58.62	1.072e-02\\
58.32	7.559e-03\\
};
\addlegendentry{RK32}

\end{axis}

\end{tikzpicture}%
  \caption{Results for the simulation of 2D DIB model~\eqref{eq:dib} with
  $N=200^2$ spatial discretization points. The number
  of time steps (or input tolerance) for each integrator is reported in
  Table~\ref{tab:DIB_ts}. The final simulation time is $T=2.5$.}
  \label{fig:DIB_cpudiag}
\end{figure}%
\begin{figure}[htb!]
  \centering
  \input{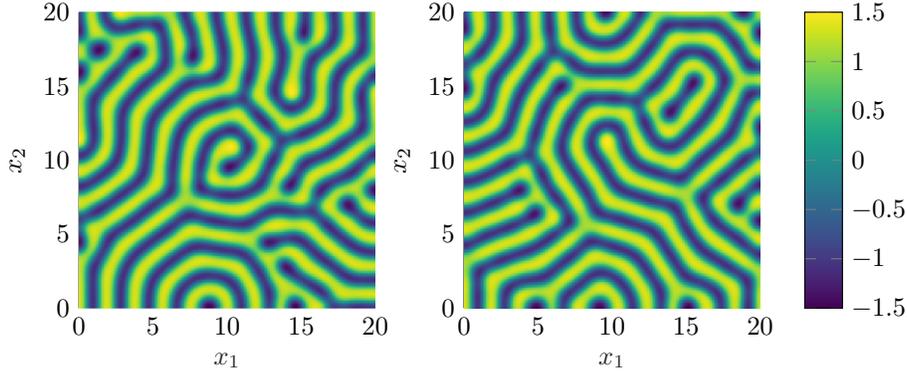}
  \caption{Turing pattern ($u$ component)for 2D DIB model~\eqref{eq:dib} obtained at final time
    $T=100$ with $N=200^2$ space discretization points with ETD2RKds and
    with time step size $\tau=2\mathrm{e}{-3}$ (left plot). The simulation
    wall-clock time is {\color{black}219} seconds. We observe
    the same pattern also by using a smaller time step
    $\tau=1.25\mathrm{e}{-3}$. However, if a larger time step
    size $\tau=2\mathrm{e}{-2}$ is used, a {\color{black} different labyrinth}
    appears (right plot).}
  \label{fig:DIB_Turing}
\end{figure}
\begin{figure}
  \centering
  \input{img/dib_mean.tex}
  \input{img/dib_increment.tex}
  \caption{{\color{black} Indicators for the time dynamics of 2D 
    DIB model~\eqref{eq:dib} solved up to final time $T=100$ 
    with  ETD2RKds
    and different time step sizes. The top plot refers to the spatial mean
    $\langle \bU_n \rangle$ while the bottom plot depicts the time increment
    $\lVert\bU_{n+1}-\bU_n \rVert_\mathrm{F}$}.}
  \label{fig:DIB_indicators}
\end{figure}
\end{subfigures}%
The plot is again divided into two parts. On the right we have the
more expensive methods, with an overall wall-clock time of roughly 100 seconds.
In particular, the ETD-RDP-IF method and the ETD2RK method do not
perform better than the explicit RK23 method. This means that for the
considered
parameters and spatial discretization the problem does not
appear to be excessively
stiff. In fact, the implicit method is roughly five times faster
than the explicit one and with performances similar to the Lawson2b method.
The ETD2RKds method is the fastest among all and, in particular, ten times
faster than
ETD2RK. {\color{black}In this experiment, the average wall-clock time of 
ETD2RKds is $4.5\mathrm{e}{-3}$ seconds per time step.}

Then, we integrate the system up to $T=100$ with ETD2RKds and 
50000 time steps. The observed steady Turing
pattern is shown in Figure~\ref{fig:DIB_Turing}, left plot.
The wall-clock time
needed for this experiment is {\color{black}219} seconds.
It was observed in Reference~\cite{DASS20} that a spatial domain
not large enough
or a number of discretization points too small
can prevent from {\color{black} clearly detecting the labyrinth Turing pattern.
  In our experiment, where the domain and the number of
  discretization points have been
properly chosen, we notice that if we take
5000 time steps instead of 50000, a clear steady labyrinth pattern still appears
(see Figure~\ref{fig:DIB_Turing}, right plot, and the indicators in
Figure~\ref{fig:DIB_indicators}).
Instead, if we increase to 80000 the number of time steps, the first
labyrinth pattern shows up again.}
\subsection{Three-dimensional
  advective Schnakenberg system}\label{sec:Schnakenberg_3d}
The diffusion--reaction Schnakenberg system
models the limit cycle behaviors of two-component
chemical reactions.
An advection term was introduced in Reference~\cite{PSM95} in order
to study its effect on patterns
(see also Reference~\cite{BKW18}). We therefore study the system
\begin{equation}\label{eq:Schnakenberg_3d}
  \left\{
  \begin{aligned}
    \partial_t u&=
    \delta^u\Delta u-\alpha^u(\partial_{x_1} u+\partial_{x_2} u+\partial_{x_3} u)+\rho(a^u-u+u^2v),\\
    \partial_t v&=
    \delta^v \Delta v-\alpha^v(\partial_{x_1} v+\partial_{x_2} v+\partial_{x_3} v)+\rho(a^v-u^2v).
  \end{aligned}\right.
\end{equation}
The computational domain is $\Omega=[0,1]^3$, the advection
parameters are $\alpha^u=\alpha^v=0.01$, the diffusion
parameters are $\delta^u=0.05$ and $\delta^v=1$,
the reaction parameter is $\rho=100$, and the
concentration parameters are $a^u=0.1305$ and $a^v=0.7695$,
respectively. As initial solution we take
\begin{equation*}
  \left\{
  \begin{aligned}
    u_0(x_1,x_2,x_3)&=a^u+a^v+10^{-5}\cdot\rme^{-100((x_1-1/3)^2+(x_2-1/2)^2+(x_3-1/3)^2)},\\
    v_0(x_1,x_2,x_3)&=\frac{a^v}{(a^u+a^v)^2},
  \end{aligned}\right.
\end{equation*}
which corresponds to a small deviation from the steady state
solution $(u_\mathrm{e},v_\mathrm{e})=(a^u+a^v,a^v/(a^u+a^v)^2)$. We integrate the system up to final
time
$T=0.4$, with the numbers of time steps (or input tolerances) given in
Table~\ref{tab:schnak_3d}. The degrees of freedom in space are $N=80^3$.
The detailed results are given in the table and depicted in Figure~\ref{fig:schnak3dadv}.
\begin{table}[htb!]
  \centering
  {\small
  \begin{tabular}{llllllll}
    \hline
    \multicolumn{4}{l}{ETD2RKds} & \multicolumn{4}{l}{ETD2RK}\\
        steps & time (s) & error & order &     steps & time (s) & error & order\\
    \hline
    50 & 1.35 & $2.34\rme{-3}$ & ---  & 50 & 10.76 & $1.48\rme{-3}$ & ---     \\
    150 & 3.84 & $3.19\rme{-4}$ & 1.81 & 150 & 25.34 & $2.79\rme{-4}$ & 1.52 \\
    250 & 6.74 & $1.26\rme{-4}$ & 1.82 & 250 & 38.71 & $1.14\rme{-4}$ & 1.76 \\
    350 & 9.70 & $6.64\rme{-5}$ & 1.91 & 350 & 51.45 & $6.04\rme{-5}$ & 1.88 \\
    \hline
    \hline
    \multicolumn{4}{l}{Lawson2b} & \multicolumn{4}{l}{ETD-RDP-IF}\\
        steps & time (s) & error & order &     steps & time (s)& error & order\\
    \hline
    400 & 6.10 & $6.90\rme{-4}$ & ---  & 200 & 52.04 & $9.98\rme{-4}$ & --- \\
    800 & 12.16 & $1.89\rme{-4}$ & 1.87 & {\color{black}450} & 116.11 & $2.59\rme{-4}$ & 1.66\\
    1200 & 18.45 & $8.59\rme{-5}$ & 1.94 & {\color{black}700} & 181.09 & $1.17\rme{-4}$ & 1.79\\
    1600 & 24.89 & $4.81\rme{-5}$ & 2.01 & {\color{black}950} & 245.61 & $6.61\rme{-5}$ & 1.88\\
    \hline
    \end{tabular}
    \begin{tabular}{lll}
    \hline
    \multicolumn{3}{l}{RK32}\\
    tolerance & time (s) & error \\
    \hline
    $8\mathrm{e}{-3}$ & 1075.21 & $8.12\rme{-3}$\\
    $4\mathrm{e}{-3}$ & 1033.12 & $3.79\rme{-3}$\\
    $8\mathrm{e}{-4}$ & 1026.94 & $8.60\rme{-4}$\\
    $4\mathrm{e}{-4}$ & 1126.25 & $2.03\rme{-4}$\\
    \hline
  \end{tabular}}%
  \caption{{\color{black}Number of time steps (or input tolerance), wall-clock time (in seconds),
  relative error at final time, and observed numerical order of convergence
  for the solution of the 3D advective
  Schnakenberg model~\eqref{eq:Schnakenberg_3d} 
      up to $T=0.4$ with different integrators.
      The DIRK23 method interrupted due to
      excessive memory requirements with tolerance $1\mathrm{e}{-1}$.
See also Figure~\ref{fig:schnak3dadv} for a graphical representation.}}
  \label{tab:schnak_3d}
\end{table}
\begin{subfigures}
\begin{figure}[htb!]
  \centering
%
%
\definecolor{mycolor1}{rgb}{1.00000,0.00000,1.00000}%
\begin{tikzpicture}

\begin{axis}[%
width=2.5in,
height=2in,
at={(0.758in,0.481in)},
scale only axis,
xmode=log,
xmin=1,
xmax=5000,
xminorticks=true,
ymode=log,
ymin=1e-05,
ymax=3e-2,
xlabel = {Wall-clock time (s)},
ylabel = {Relative error},
yminorticks=true,
axis background/.style={fill=white},
legend style={at={(0.02,0.725)}, anchor=south west, legend cell align=left, align=left, draw=white!15!black, font=\scriptsize},
legend columns = 2
]
\addplot [color=mycolor1,line width=1.5pt, mark size = 3pt, mark=+, mark options={solid, mycolor1}]
  table[row sep=crcr]{%
1.349538	0.00233518293850402\\
3.837403	0.000319349860140726\\
6.735411	0.000126113102253263\\
9.700735	6.63960004225219e-05\\
};
\addlegendentry{ETD2RKds}

\addplot [color=blue, line width=1.5pt, mark size = 3pt,mark=o, mark options={solid, blue}]
  table[row sep=crcr]{%
10.759346	0.00148017947603509\\
25.335923	0.000279184568193472\\
38.713644	0.000113667748736843\\
51.451994	6.04337009884737e-05\\
};
\addlegendentry{ETD2RK}

\addplot [color=red,line width=1.5pt, mark size = 3pt, mark=triangle, mark options={solid, rotate=270, red}]
  table[row sep=crcr]{%
6.101923	0.000690407285509494\\
12.15621	0.000188954464550619\\
18.451787	8.5888841287558e-05\\
24.888401	4.81198193585188e-05\\
};
\addlegendentry{Lawson2b}


\addplot [color=green, line width=1.5pt, mark size = 3pt, mark=asterisk, mark options={solid, green}]
  table[row sep=crcr]{%
52.04	9.9758e-04\\%
116.11	2.590e-04\\%
181.09	1.172e-04\\%
245.61	6.605e-05\\%
};
\addlegendentry{ETD-RDP-IF}

\addplot [color=black,line width=1.5pt, mark size = 3pt, mark=triangle, mark options={solid, rotate=90, black}]
  table[row sep=crcr]{%
1.0752e+03	0.0081\\
1.0331e+03	3.7883e-03\\
1.0269e+03	8.6011e-04\\
1.1262e+03	2.0307e-04\\
};
\addlegendentry{RK32}

\end{axis}

\end{tikzpicture}%
  \caption{Results for the simulation of the 3D advective Schnakenberg
  model~\eqref{eq:Schnakenberg_3d} with
  $N=80^3$ spatial discretization points. The number
  of time steps (or input tolerance) for each integrator is reported in
  Table~\ref{tab:schnak_3d}. The final simulation time is $T=0.4$.}
  \label{fig:schnak3dadv}
\end{figure}
\begin{figure}[htb!]
  \centering
  \input{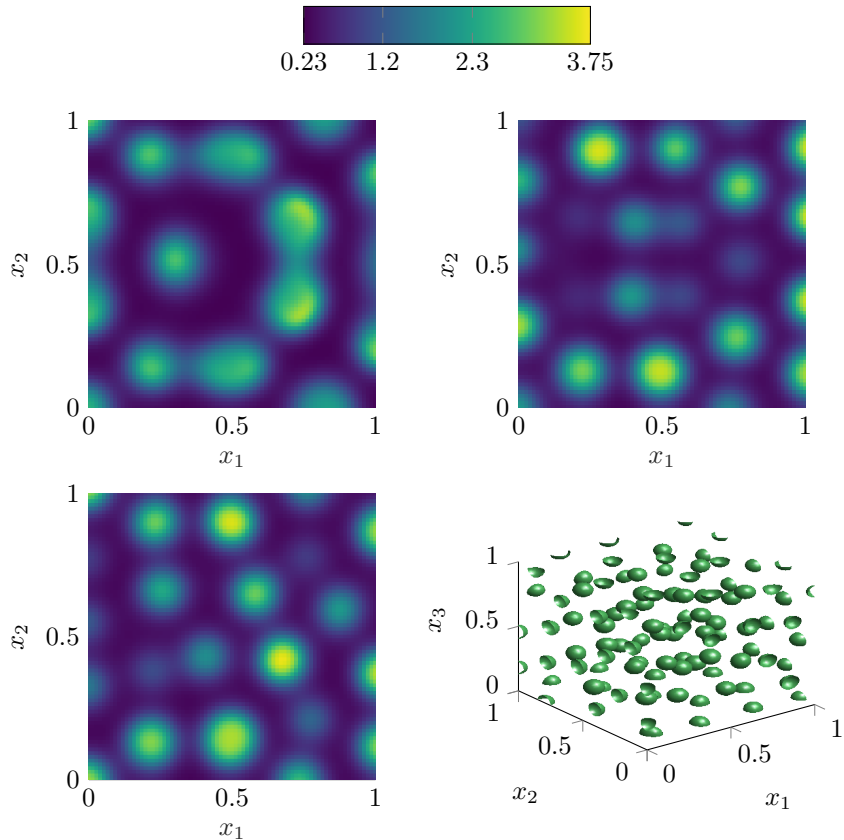}
  \caption{Spot-like
    pattern ($u$ component) for the 3D
    advective Schnakenberg model~\eqref{eq:Schnakenberg_3d}
    {\color{black}with $N=80^3$ space discretization points and ETD2RKds at final
    time $T=0.8$ (top left), $T=8$ (top right) and $T=80$ (bottom left).
    The simulation wall-clock time is 3, 28 and 275 seconds, respectively.
    The time step size employed is $\tau=8\mathrm{e}{-3}$, and 
    the slices correspond to $x_3=0.49$.
    The common colorbar is displayed at the top. We report also the 
    isosurface of level $2.8$ at $T=80$ (bottom right plot).}}
  \label{fig:schnak3dpatt}
\end{figure}
\end{subfigures}

As in the previous three-dimensional example, the implicit Runge--Kutta method
DIRK23 could
not terminate due to excessive memory requirements.
Concerning the other methods, we 
{\color{black}observe that explicit Runge--Kutta scheme performs poorly in 
terms of wall-clock time. Better performances can be gradually obtained by
employing the ETD-RDP-IF method, the ETD2RK scheme and the Lawson2b integrator.}
Overall, the best method turns out to be again ETD2RKds, {\color{black} with an 
average wall-clock time per time step of $2.7\mathrm{e}{-2}$ seconds}.

Then, we repeat the experiment up to final {\color{black}times $T=0.8$, $T=8$, and
$T=80$} with again $N=80^3$
degrees of freedom in space and using as time integrator ETD2RKds with
$\tau=8\mathrm{e}{-3}$ {\color{black} (i.e., performing 100, 1000 and 10000 time steps,
respectively)}.
As already observed in the literature (see Reference~\cite{BKW18}),
the initial condition evolves to
a spot-like pattern, {\color{black}as}
depicted in Figure~\ref{fig:schnak3dpatt}.
The overall simulation {\color{black}times are roughly 3, 28, and 275 seconds
for the final times $T=0.8$, $T=8$, and $T=80$, respectively}.

\subsection{Three-dimensional
  advective Brusselator system}\label{sec:Brusselator}
We consider the advective Brusselator system
\begin{equation}\label{eq:Brusselator}
  \left\{
  \begin{aligned}
    \partial_t u&=
    \delta^u\Delta u-\alpha^u(\partial_{x_1}u + \partial_{x_2}u+\partial_{x_3}u)+
    u^2v-(a_1^u+1)u+a_2^u,\\
    \partial_t v&=
    \delta^v \Delta v-\alpha^v(\partial_{x_1}v + \partial_{x_2}v+\partial_{x_3}v)
    -u^2v+a_1^uu,
  \end{aligned}\right.
\end{equation}
presented
in Reference~\cite{BKW18}. The unknowns $u$ and $v$ represent
the concentration of the activator and inhibitor in the
chemical reaction, respectively.
The spatial domain is $\Omega=[0,1]^3$,
the diffusion coefficients are $\delta^u=0.01$, $\delta^v=0.02$,
the advection coefficients are $a^u_1=1$,
and the remaining parameters are
$a_2^u=2$, $\alpha^u = \alpha^v = 0.1$. The space discretization is
performed using $N=64^3$ points, and
the initial datum is 
$u_0(x_1,x_2,x_3) = 1 + \sin(2\pi x_1)\sin(2\pi x_2)\sin(2\pi x_3)$
and $v_0(x_1,x_2,x_3)=3$. The number of steps for each integrator is reported in
Table~\ref{tab:bruss}, as well as the detailed outcome of the experiment 
for final time $T=1$. Similarly to the previous examples, we also present the
precision diagram in Figure~\ref{fig:bruss}.
\begin{table}[htb!]
  \centering
  {\small
  \begin{tabular}{llllllll}
    \hline
    \multicolumn{4}{l}{ETD2RKds} & \multicolumn{4}{l}{ETD2RK}\\
        steps & time (s) & error & order &     steps & time (s) & error & order\\
    \hline
    50 & 0.48 & $3.46\rme{-4}$ & ---  & 50 & 4.76 & $3.46\rme{-4}$ & ---     \\
    100 & 0.94 & $7.94\rme{-5}$ & 2.13 & 100 & 7.84 & $7.93\rme{-5}$ & 2.12 \\
    150 & 1.25 & $3.43\rme{-5}$ & 2.07 & 150 & 11.38 & $3.43\rme{-5}$ & 2.07 \\
    200 & 1.82 & $1.90\rme{-5}$ & 2.05 & 200 & 13.62 & $1.90\rme{-5}$ & 2.05 \\
    \hline
    \hline
    \multicolumn{4}{l}{Lawson2b} & \multicolumn{4}{l}{ETD-RDP-IF}\\
        steps & time (s) & error & order &     steps & time (s)& error & order\\
    \hline
    50 & 0.32 & $3.43\rme{-4}$ & ---  & 50 & 6.59 & $3.35\rme{-4}$ & --- \\
    100 & 0.63 & $7.85\rme{-5}$ & 2.13 & 100 & 12.91 & $7.68\rme{-5}$ & 2.12\\
    150 & 0.81 & $3.40\rme{-5}$ & 2.07 & 150 & 19.24 & $3.32\rme{-5}$ & 2.07\\
    200 & 1.31 & $1.88\rme{-5}$ & 2.05 & 200 & 25.39 & $1.84\rme{-5}$ & 2.05\\
    \hline
    \end{tabular}
    \begin{tabular}{lll}
    \hline
    \multicolumn{3}{l}{RK32}\\
    tolerance & time (s) & error \\
    \hline
    $8\mathrm{e}{-4}$ & 16.08 & $8.73\rme{-4}$\\
    $6\mathrm{e}{-4}$ & 16.65 & $1.11\rme{-4}$\\
    $1\mathrm{e}{-4}$ & 16.54 & $5.48\rme{-5}$\\
    $5\mathrm{e}{-5}$ & 16.64 & $2.43\rme{-5}$\\
    \hline
  \end{tabular}}%
  \caption{{\color{black}Number of time steps (or input tolerance), wall-clock time (in seconds),
  relative error at final time, and observed numerical order of convergence
      for the solution of the 3D advective Brusselator model~\eqref{eq:Brusselator} 
      up to $T=1$ with different integrators.
      The DIRK23 method interrupted due to
      excessive memory requirements with tolerance $1\mathrm{e}{-1}$.
See also Figure~\ref{fig:bruss} for a graphical representation.}}
  \label{tab:bruss}
\end{table}
\begin{subfigures}
\begin{figure}[htb!]
  \centering
%
%
\definecolor{mycolor1}{rgb}{1.00000,0.00000,1.00000}%
\begin{tikzpicture}

\begin{axis}[%
width=2.5in,
height=2in,
at={(0.769in,0.477in)},
scale only axis,
xmode=log,
xmin=0.1,
xmax=50,
xminorticks=true,
ymode=log,
ymin=1e-05,
ymax=0.002,
xlabel = {Wall-clock time (s)},
ylabel = {Relative error},
yminorticks=true,
axis background/.style={fill=white},
legend style={at={(0.025,0.725)}, anchor=south west, legend cell align=left, align=left, draw=white!15!black, font=\scriptsize},
legend columns=2
]
\addplot [color=mycolor1,line width=1.5pt, mark size = 3pt, mark=+, mark options={solid, mycolor1}]
  table[row sep=crcr]{%
0.477654	0.000346282011062941\\
0.937761	7.93663151979711e-05\\
1.24909	3.43302300081084e-05\\
1.815443	1.90497910690793e-05\\
};
\addlegendentry{ETD2RKds}

\addplot [color=blue, line width=1.5pt, mark size = 3pt,mark=o, mark options={solid, blue}]
  table[row sep=crcr]{%
4.758315	0.000346075110471962\\
7.84456300000001	7.93195240999125e-05\\
11.375277	3.4310181077751e-05\\
13.622031	1.9038724610664e-05\\
};
\addlegendentry{ETD2RK}

\addplot [color=red,line width=1.5pt, mark size = 3pt, mark=triangle, mark options={solid, rotate=270, red}]
  table[row sep=crcr]{%
0.317235	0.000342788309638027\\
0.626430000000001	7.85387500125735e-05\\
0.808573	3.39671219130416e-05\\
1.310796	1.88466882792953e-05\\
};
\addlegendentry{Lawson2b}


\addplot [color=green, line width=1.5pt, mark size = 3pt, mark=asterisk, mark options={solid, green}]
  table[row sep=crcr]{%
6.59	0.000335324335914824\\
12.91	7.678114503436e-05\\
19.24	3.32016964284683e-05\\
25.39	1.84205614353252e-05\\
};
\addlegendentry{ETD-RDP-IF}

\addplot [color=black, line width=1.5pt, mark size = 3pt, mark=triangle, mark options={solid, rotate=90, black}]
  table[row sep=crcr]{%
16.083921	0.00087282815749488\\
16.645004	0.000110812476916986\\
16.540105	5.47544899262375e-05\\
16.642905	2.43097258338379e-05\\
};
\addlegendentry{RK32}

\end{axis}

\end{tikzpicture}%
  \caption{Results for the simulation of the 3D advective Brusselator
  model~\eqref{eq:Brusselator} with
  $N=64^3$ spatial discretization points. The number
  of time steps (or input tolerance) for each integrator is reported in
  Table~\ref{tab:bruss}. The final simulation time
  is $T=1$.}
  \label{fig:bruss}
\end{figure}
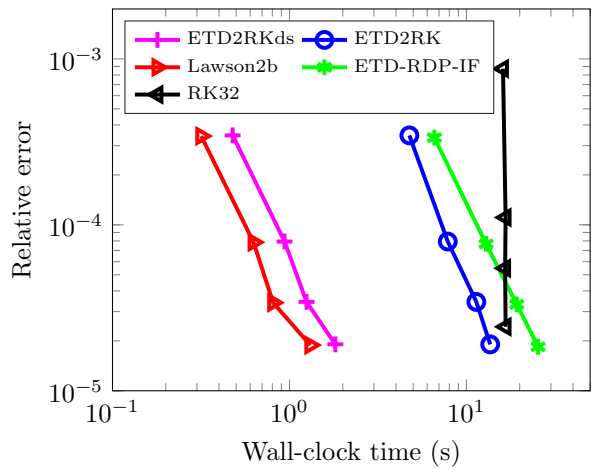
\begin{figure}[htb!]
  \centering
  \input{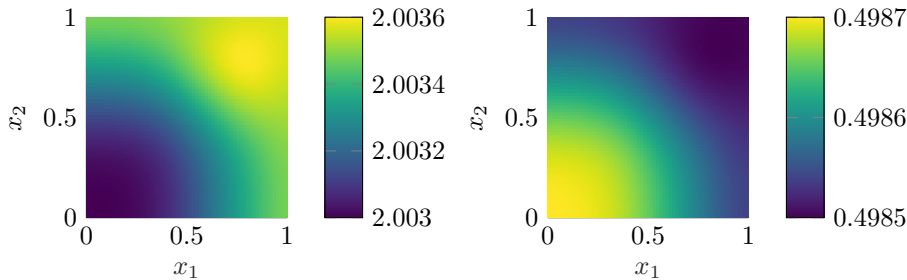}
  \caption{Equilibrium state of the 3D advective Brusselator
    model~\eqref{eq:Brusselator}
    obtained at final time
    $T=5$ with $N=64^3$ space discretization points and ETD2RKds. The
    time step size employed is $\tau=5\mathrm{e}{-2}$. The simulation wall-clock
    time is 1 second.
    The reported slices correspond to $x_3=1$
    for the component $u$
    (left plot) and the component $v$ (right plot).}
  \label{fig:Brusselator_patt}
\end{figure}
\end{subfigures}

Again, the plot is split into two parts. Concerning the most efficient schemes,
in this experiment the Lawson2b method is slightly faster than the ETD2RKds 
scheme, while reaching the same set of accuracies.
{\color{black}The average wall-clock time per time step of
ETD2RKds is $9.1\mathrm{e}{-3}$ seconds.}
All the remaining methods are almost one order of magnitude slower, with the
DIRK23 method not able to produce an approximation due to excessive
memory requirements.

We then simulate the system
up to $T=5$ with 100 time steps with ETD2RKds
(the overall computational time is about 1 second).  As already observed
in Reference~\cite{BKW18}, the solution approaches the equilibrium
state $(u_\mathrm{e},v_\mathrm{e})=(a_2^u,a_1^u/a_2^u)$, see
Figure~\ref{fig:Brusselator_patt}.

\section{Conclusions}\label{sec:conc}
In this paper, we show how it is possible to effectively exploit the Kronecker
sum structure {\color{black}for the time integration of semidiscretized}
two-component systems of coupled advection--diffusion--reaction equations.
The proposed second order
exponential-type
time marching scheme, which is based on a directional splitting of the involved
matrix functions and named ETD2RKds, is shown to outperform well-established
techniques on a variety
of physically relevant models from the literature, such as two-dimensional
Schnakenberg, FitzHugh--Nagumo, DIB, 
and three-dimensional FitzHugh--Nagumo, advective
Schnakenberg, and
advective Brusselator models. The procedure {\color{black}is able to capture
the formation of Turing patterns, easily extends to models with any
number of components,
and can be applied in more than three spatial dimensions.}
As future work, we plan to
investigate and
perform simulations at HPC level (i.e., using server multi-core CPUs and GPUs)
to further enhance the performances of the directional splitting procedure.
\bibliography{CC24}

\begin{thebibliography}{39}
\expandafter\ifx\csname natexlab\endcsname\relax\def\natexlab#1{#1}\fi
\providecommand{\url}[1]{\texttt{#1}}
\providecommand{\href}[2]{#2}
\providecommand{\path}[1]{#1}
\providecommand{\DOIprefix}{doi:}
\providecommand{\ArXivprefix}{arXiv:}
\providecommand{\URLprefix}{URL: }
\providecommand{\Pubmedprefix}{pmid:}
\providecommand{\doi}[1]{\href{http://dx.doi.org/#1}{\path{#1}}}
\providecommand{\Pubmed}[1]{\href{pmid:#1}{\path{#1}}}
\providecommand{\bibinfo}[2]{#2}
\ifx\xfnm\relax \def\xfnm[#1]{\unskip,\space#1}\fi
\bibitem[{Al-Mohy and Higham(2010)}]{AMH10}
\bibinfo{author}{Al-Mohy, A.H.}, \bibinfo{author}{Higham, N.J.},
  \bibinfo{year}{2010}.
\newblock \bibinfo{title}{A new scaling and squaring algorithm for the matrix
  exponential}.
\newblock \bibinfo{journal}{SIAM J. Matrix Anal. Appl.} \bibinfo{volume}{31},
  \bibinfo{pages}{970--989}.
\bibitem[{Alla et~al.(2023)Alla, Monti and Sgura}]{AMS23}
\bibinfo{author}{Alla, A.}, \bibinfo{author}{Monti, A.},
  \bibinfo{author}{Sgura, I.}, \bibinfo{year}{2023}.
\newblock \bibinfo{title}{Adaptive {POD}-{DEIM} correction for {T}uring pattern
  approximation in reaction--diffusion {PDE} systems}.
\newblock \bibinfo{journal}{J. Numer. Math.} \bibinfo{volume}{31},
  \bibinfo{pages}{205--229}.
\bibitem[{Alonso et~al.(2023)Alonso, Ib\'a{\~n}ez, Defez and
  Alonso-Jord\'a}]{AIDAJ23}
\bibinfo{author}{Alonso, J.M.}, \bibinfo{author}{Ib\'a{\~n}ez, J.},
  \bibinfo{author}{Defez, E.}, \bibinfo{author}{Alonso-Jord\'a, P.},
  \bibinfo{year}{2023}.
\newblock \bibinfo{title}{Euler polynomials for the matrix exponential
  approximation}.
\newblock \bibinfo{journal}{J. Comput. Appl. Math.} \bibinfo{volume}{425},
  \bibinfo{pages}{115074}.
\bibitem[{Asante-Asamani et~al.(2020)Asante-Asamani, Kleefeld and
  Wade}]{AAKW20}
\bibinfo{author}{Asante-Asamani, E.O.}, \bibinfo{author}{Kleefeld, A.},
  \bibinfo{author}{Wade, B.A.}, \bibinfo{year}{2020}.
\newblock \bibinfo{title}{A second-order exponential time differencing scheme
  for non-linear reaction-diffusion systems with dimensional splitting}.
\newblock \bibinfo{journal}{J. Comput. Phys.} \bibinfo{volume}{415},
  \bibinfo{pages}{109490}.
\bibitem[{Ben~Tahar et~al.(2023)Ben~Tahar, Mu{\~n}oz, Shefelbine and
  Comellas}]{BTMSC23}
\bibinfo{author}{Ben~Tahar, S.}, \bibinfo{author}{Mu{\~n}oz, J.J.},
  \bibinfo{author}{Shefelbine, S.J.}, \bibinfo{author}{Comellas, E.},
  \bibinfo{year}{2023}.
\newblock \bibinfo{title}{Turing pattern prediction in three-dimensional
  domains: the role of initial conditions and growth}.
\newblock \bibinfo{journal}{bioRxiv 2023.03.29.534782} .
\bibitem[{Berland et~al.(2005)Berland, Skaflestad and Wright}]{BSW05}
\bibinfo{author}{Berland, H.}, \bibinfo{author}{Skaflestad, B.},
  \bibinfo{author}{Wright, W.M.}, \bibinfo{year}{2005}.
\newblock \bibinfo{title}{{EXPINT} --- {A} {MATLAB} package for exponential
  integrators}.
\newblock \bibinfo{type}{Technical Report} \bibinfo{number}{4}. Norwegian
  University of Science and Technology.
\bibitem[{Bhatt et~al.(2018)Bhatt, Khaliq and Wade}]{BKW18}
\bibinfo{author}{Bhatt, H.P.}, \bibinfo{author}{Khaliq, A.Q.M.},
  \bibinfo{author}{Wade, B.A.}, \bibinfo{year}{2018}.
\newblock \bibinfo{title}{Efficient {K}rylov-based exponential time
  differencing method in application to {3D} advection-diffusion-reaction
  systems}.
\newblock \bibinfo{journal}{Appl. Math. Comput.} \bibinfo{volume}{338},
  \bibinfo{pages}{260--273}.
\bibitem[{Bogacki and Shampine(1989)}]{BS89}
\bibinfo{author}{Bogacki, P.}, \bibinfo{author}{Shampine, L.F.},
  \bibinfo{year}{1989}.
\newblock \bibinfo{title}{A 3(2) pair of {R}unge -- {K}utta formulas}.
\newblock \bibinfo{journal}{Appl. Math. Lett.} \bibinfo{volume}{2},
  \bibinfo{pages}{321--325}.
\bibitem[{Bozzini et~al.(2013)Bozzini, Lacitignola and Sgura}]{BLS13}
\bibinfo{author}{Bozzini, B.}, \bibinfo{author}{Lacitignola, D.},
  \bibinfo{author}{Sgura, I.}, \bibinfo{year}{2013}.
\newblock \bibinfo{title}{Spatio-temporal organization in alloy
  electrodeposition: a morphochemical mathematical model and its experimental
  validation}.
\newblock \bibinfo{journal}{J. Solid State Electrochem.} \bibinfo{volume}{17},
  \bibinfo{pages}{467--479}.
\bibitem[{Caliari and Cassini(2023)}]{CC23}
\bibinfo{author}{Caliari, M.}, \bibinfo{author}{Cassini, F.},
  \bibinfo{year}{2023}.
\newblock \bibinfo{title}{Direction splitting of $\varphi$-functions in
  exponential integrators for $d$-dimensional problems in {K}ronecker form}.
\newblock \bibinfo{journal}{J. Approx. Softw.} \bibinfo{note}{Accepted for
  publication}.
\bibitem[{Caliari et~al.(2022)Caliari, Cassini, Einkemmer, Ostermann and
  Zivcovich}]{CCEOZ22}
\bibinfo{author}{Caliari, M.}, \bibinfo{author}{Cassini, F.},
  \bibinfo{author}{Einkemmer, L.}, \bibinfo{author}{Ostermann, A.},
  \bibinfo{author}{Zivcovich, F.}, \bibinfo{year}{2022}.
\newblock \bibinfo{title}{A $\mu$-mode integrator for solving evolution
  equations in {K}ronecker form}.
\newblock \bibinfo{journal}{J. Comput. Phys.} \bibinfo{volume}{455},
  \bibinfo{pages}{110989}.
\bibitem[{Caliari et~al.(2023a)Caliari, Cassini and Zivcovich}]{CCZ23}
\bibinfo{author}{Caliari, M.}, \bibinfo{author}{Cassini, F.},
  \bibinfo{author}{Zivcovich, F.}, \bibinfo{year}{2023}a.
\newblock \bibinfo{title}{{BAMPHI}: {M}atrix-free and transpose-free action of
  linear combinations of $\varphi$-functions from exponential integrators}.
\newblock \bibinfo{journal}{J. Comput. Appl. Math.} \bibinfo{volume}{423},
  \bibinfo{pages}{114973}.
\bibitem[{Caliari et~al.(2023b)Caliari, Cassini and Zivcovich}]{CCZ23phi}
\bibinfo{author}{Caliari, M.}, \bibinfo{author}{Cassini, F.},
  \bibinfo{author}{Zivcovich, F.}, \bibinfo{year}{2023}b.
\newblock \bibinfo{title}{A $\mu$-mode approach for exponential integrators:
  actions of $\varphi$-functions of {K}ronecker sums}.
\newblock \bibinfo{journal}{arXiv:2210.07667 [math.NA]} .
\bibitem[{Caliari et~al.(2023c)Caliari, Cassini and Zivcovich}]{CCZ23kp}
\bibinfo{author}{Caliari, M.}, \bibinfo{author}{Cassini, F.},
  \bibinfo{author}{Zivcovich, F.}, \bibinfo{year}{2023}c.
\newblock \bibinfo{title}{A $\mu$-mode {BLAS} approach for multidimensional
  tensor-structured problems}.
\newblock \bibinfo{journal}{Numer. Algorithms} \bibinfo{volume}{92},
  \bibinfo{pages}{2483--2508}.
\bibitem[{Caliari and Zivcovich(2019)}]{CZ19}
\bibinfo{author}{Caliari, M.}, \bibinfo{author}{Zivcovich, F.},
  \bibinfo{year}{2019}.
\newblock \bibinfo{title}{On-the-fly backward error estimate for matrix
  exponential approximation by {T}aylor algorithm}.
\newblock \bibinfo{journal}{J. Comput. Appl. Math.} \bibinfo{volume}{346},
  \bibinfo{pages}{532--548}.
\bibitem[{Croci and Mu{\~n}oz-Matute(2023)}]{CMM23}
\bibinfo{author}{Croci, M.}, \bibinfo{author}{Mu{\~n}oz-Matute, J.},
  \bibinfo{year}{2023}.
\newblock \bibinfo{title}{Exploiting {K}ronecker structure in exponential
  integrators: {F}ast approximation of the action of $\varphi$-functions of
  matrices via quadrature}.
\newblock \bibinfo{journal}{J. Comput. Sci.} \bibinfo{volume}{67},
  \bibinfo{pages}{101966}.
\bibitem[{D’Autilia et~al.(2020)D’Autilia, Sgura and Simoncini}]{DASS20}
\bibinfo{author}{D’Autilia, M.C.}, \bibinfo{author}{Sgura, I.},
  \bibinfo{author}{Simoncini, V.}, \bibinfo{year}{2020}.
\newblock \bibinfo{title}{Matrix-oriented discretization methods for
  reaction--diffusion {PDE}s: {C}omparisons and applications}.
\newblock \bibinfo{journal}{Comput. Math. with Appl.} \bibinfo{volume}{79},
  \bibinfo{pages}{2067--2085}.
\bibitem[{Gambino et~al.(2019)Gambino, Lombardo, Rubino and
  Sammartino}]{GLRS19}
\bibinfo{author}{Gambino, G.}, \bibinfo{author}{Lombardo, M.C.},
  \bibinfo{author}{Rubino, G.}, \bibinfo{author}{Sammartino, M.},
  \bibinfo{year}{2019}.
\newblock \bibinfo{title}{Pattern selection in the 2{D} {F}itz{H}ugh--{N}agumo
  model}.
\newblock \bibinfo{journal}{Ric. di Mat.} \bibinfo{volume}{68},
  \bibinfo{pages}{535--549}.
\bibitem[{Gaudreault et~al.(2018)Gaudreault, Rainwater and Tokman}]{GRT18}
\bibinfo{author}{Gaudreault, S.}, \bibinfo{author}{Rainwater, G.},
  \bibinfo{author}{Tokman, M.}, \bibinfo{year}{2018}.
\newblock \bibinfo{title}{{KIOPS: A} fast adaptive {K}rylov subspace solver for
  exponential integrators}.
\newblock \bibinfo{journal}{J. Comput. Phys.} \bibinfo{volume}{372},
  \bibinfo{pages}{236--255}.
\bibitem[{Hairer and Wanner(1996)}]{HW96}
\bibinfo{author}{Hairer, E.}, \bibinfo{author}{Wanner, G.},
  \bibinfo{year}{1996}.
\newblock \bibinfo{title}{Solving Ordinary Differential Equations II: Stiff and
  Differential-Algebraic Problems}. volume~\bibinfo{volume}{14} of
  \textit{\bibinfo{series}{Springer Series in Computational Mathematics}}.
\newblock \bibinfo{edition}{second} ed., \bibinfo{publisher}{Springer Berlin}.
\bibitem[{Hochbruck et~al.(2020)Hochbruck, Leibold and Ostermann}]{HLO20}
\bibinfo{author}{Hochbruck, M.}, \bibinfo{author}{Leibold, J.},
  \bibinfo{author}{Ostermann, A.}, \bibinfo{year}{2020}.
\newblock \bibinfo{title}{On the convergence of {L}awson methods for semilinear
  stiff problems}.
\newblock \bibinfo{journal}{Numer. Math.} \bibinfo{volume}{145},
  \bibinfo{pages}{553--580}.
\bibitem[{Hochbruck and Ostermann(2010)}]{HO10}
\bibinfo{author}{Hochbruck, M.}, \bibinfo{author}{Ostermann, A.},
  \bibinfo{year}{2010}.
\newblock \bibinfo{title}{Exponential integrators}.
\newblock \bibinfo{journal}{Acta Numer.} \bibinfo{volume}{19},
  \bibinfo{pages}{209--286}.
\bibitem[{Hosea and Shampine(1996)}]{HS96}
\bibinfo{author}{Hosea, M.E.}, \bibinfo{author}{Shampine, L.F.},
  \bibinfo{year}{1996}.
\newblock \bibinfo{title}{Analysis and implementation of {TR-BDF2}}.
\newblock \bibinfo{journal}{Appl. Numer. Math.} \bibinfo{volume}{20},
  \bibinfo{pages}{21--37}.
\bibitem[{Jiang and Zhang(2016)}]{JZ16}
\bibinfo{author}{Jiang, T.}, \bibinfo{author}{Zhang, Y.T.},
  \bibinfo{year}{2016}.
\newblock \bibinfo{title}{Krylov single-step implicit integration factor {WENO}
  methods for advection--diffusion--reaction equations}.
\newblock \bibinfo{journal}{J. Comput. Phys.} \bibinfo{volume}{311},
  \bibinfo{pages}{22--44}.
\bibitem[{Lawson(1967)}]{L67}
\bibinfo{author}{Lawson, J.D.}, \bibinfo{year}{1967}.
\newblock \bibinfo{title}{Generalized {R}unge-{K}utta processes for stable
  systems with large {L}ipschitz constants}.
\newblock \bibinfo{journal}{SIAM J. Numer. Anal.} \bibinfo{volume}{4},
  \bibinfo{pages}{372--380}.
\bibitem[{Li et~al.(2022)Li, Yang and Lan}]{LYL22}
\bibinfo{author}{Li, D.}, \bibinfo{author}{Yang, S.}, \bibinfo{author}{Lan,
  J.}, \bibinfo{year}{2022}.
\newblock \bibinfo{title}{Efficient and accurate computation for the
  $\varphi$‑functions arising from exponential integrators}.
\newblock \bibinfo{journal}{Calcolo} \bibinfo{volume}{59},
  \bibinfo{pages}{1--24}.
\bibitem[{Luan et~al.(2019)Luan, Pudykiewicz and Reynolds}]{LPR19}
\bibinfo{author}{Luan, V.T.}, \bibinfo{author}{Pudykiewicz, J.A.},
  \bibinfo{author}{Reynolds, D.R.}, \bibinfo{year}{2019}.
\newblock \bibinfo{title}{Further development of efficient and accurate time
  integration schemes for meteorological models}.
\newblock \bibinfo{journal}{J. Comput. Phys.} \bibinfo{volume}{376},
  \bibinfo{pages}{817--837}.
\bibitem[{Madzvamuse et~al.(2003)Madzvamuse, Wathen and Maini}]{MWM03}
\bibinfo{author}{Madzvamuse, A.}, \bibinfo{author}{Wathen, A.J.},
  \bibinfo{author}{Maini, P.K.}, \bibinfo{year}{2003}.
\newblock \bibinfo{title}{A moving grid finite element method applied to a
  model biological pattern generator}.
\newblock \bibinfo{journal}{J. Comput. Phys.} \bibinfo{volume}{190},
  \bibinfo{pages}{478--500}.
\bibitem[{Malchow et~al.(2008)Malchow, Petrovskii and Venturino}]{MPV08}
\bibinfo{author}{Malchow, H.}, \bibinfo{author}{Petrovskii, S.V.},
  \bibinfo{author}{Venturino, E.}, \bibinfo{year}{2008}.
\newblock \bibinfo{title}{Spatiotemporal {P}atterns in {E}cology and
  {E}pidemiology: {T}heory, {M}odels, and {S}imulation}.
\newblock CRC Mathematical Biology Series. \bibinfo{edition}{first} ed.,
  \bibinfo{publisher}{Chapman \& Hall}.
\bibitem[{M\"uller(2022)}]{M22}
\bibinfo{author}{M\"uller, B.}, \bibinfo{year}{2022}.
\newblock \bibinfo{title}{Investigation of {E}xponential {T}ime {D}ifferencing
  schemes for advection-diffusion-reaction problems in the presence of
  significant advection}.
\newblock Master's thesis. University of Louisiana at Lafayette.
\bibitem[{Mu{\~n}oz-Matute et~al.(2022)Mu{\~n}oz-Matute, Pardo and
  Calo}]{MMPC22}
\bibinfo{author}{Mu{\~n}oz-Matute, J.}, \bibinfo{author}{Pardo, D.},
  \bibinfo{author}{Calo, V.M.}, \bibinfo{year}{2022}.
\newblock \bibinfo{title}{Exploiting the {K}ronecker product structure of
  $\varphi$--functions in exponential integrators}.
\newblock \bibinfo{journal}{Int. J. Numer. Methods Eng.} \bibinfo{volume}{123},
  \bibinfo{pages}{2142--2161}.
\bibitem[{Neudecker(1969)}]{N69}
\bibinfo{author}{Neudecker, H.}, \bibinfo{year}{1969}.
\newblock \bibinfo{title}{A note on {K}ronecker matrix products and matrix
  equation systems}.
\newblock \bibinfo{journal}{SIAM J. Appl. Math.} \bibinfo{volume}{17},
  \bibinfo{pages}{603--606}.
\bibitem[{Perumpanani et~al.(1995)Perumpanani, Sherratt and Maini}]{PSM95}
\bibinfo{author}{Perumpanani, A.J.}, \bibinfo{author}{Sherratt, J.A.},
  \bibinfo{author}{Maini, P.K.}, \bibinfo{year}{1995}.
\newblock \bibinfo{title}{Phase differences in reaction--diffusion--advection
  systems and applications to morphogenesis}.
\newblock \bibinfo{journal}{IMA J. Appl. Math.} \bibinfo{volume}{55},
  \bibinfo{pages}{19--33}.
\bibitem[{Sastre et~al.(2019)Sastre, Ib{\'a\~n}ez and Defez}]{SID19}
\bibinfo{author}{Sastre, J.}, \bibinfo{author}{Ib{\'a\~n}ez, J.},
  \bibinfo{author}{Defez, E.}, \bibinfo{year}{2019}.
\newblock \bibinfo{title}{Boosting the computation of the matrix exponential}.
\newblock \bibinfo{journal}{Appl. Math. Comput.} \bibinfo{volume}{340},
  \bibinfo{pages}{206--220}.
\bibitem[{Schnakenberg(1979)}]{S79}
\bibinfo{author}{Schnakenberg, J.}, \bibinfo{year}{1979}.
\newblock \bibinfo{title}{Simple chemical reaction systems with limit cycle
  behaviour}.
\newblock \bibinfo{journal}{J. Theor. Biol.} \bibinfo{volume}{81},
  \bibinfo{pages}{389--400}.
\bibitem[{Sherratt and Chaplain(2001)}]{CS01}
\bibinfo{author}{Sherratt, J.A.}, \bibinfo{author}{Chaplain, M.A.J.},
  \bibinfo{year}{2001}.
\newblock \bibinfo{title}{A new mathematical model for avascular tumour
  growth}.
\newblock \bibinfo{journal}{J. Math. Biol.} \bibinfo{volume}{43},
  \bibinfo{pages}{291--312}.
\bibitem[{Singh et~al.(2023)Singh, Mittal, Thottoli and Tamsir}]{SMTT23}
\bibinfo{author}{Singh, S.}, \bibinfo{author}{Mittal, R.C.},
  \bibinfo{author}{Thottoli, S.R.}, \bibinfo{author}{Tamsir, M.},
  \bibinfo{year}{2023}.
\newblock \bibinfo{title}{High-fidelity simulations for {T}uring pattern
  formation in multi-dimensional {G}ray--{S}cott reaction-diffusion system}.
\newblock \bibinfo{journal}{Appl. Math. Comput.} \bibinfo{volume}{452},
  \bibinfo{pages}{128079}.
\bibitem[{Skaflestad and Wright(2009)}]{SW09}
\bibinfo{author}{Skaflestad, B.}, \bibinfo{author}{Wright, W.M.},
  \bibinfo{year}{2009}.
\newblock \bibinfo{title}{The scaling and modified squaring method for matrix
  functions related to the exponential}.
\newblock \bibinfo{journal}{Appl. Numer. Math.} \bibinfo{volume}{59},
  \bibinfo{pages}{783--799}.
\bibitem[{Twizell et~al.(1999)Twizell, Gumel and Cao}]{TGC99}
\bibinfo{author}{Twizell, E.H.}, \bibinfo{author}{Gumel, A.B.},
  \bibinfo{author}{Cao, Q.}, \bibinfo{year}{1999}.
\newblock \bibinfo{title}{A second-order scheme for the ``{B}russelator''
  reaction--diffusion system}.
\newblock \bibinfo{journal}{J. Math. Chem.} \bibinfo{volume}{26},
  \bibinfo{pages}{297--316}.

\end{thebibliography}
\bibliographystyle{elsarticle-harv}
\end{document}